\pdfoutput=1
\RequirePackage{ifpdf}
\ifpdf 
\documentclass[pdftex]{sigma}
\else
\documentclass{sigma}
\fi

\usepackage{tikz}
\usepackage[all]{xy}
\usepackage{nicefrac,enumitem}

\numberwithin{equation}{section}

\newtheorem{Theorem}{Theorem}[section]
\newtheorem{Corollary}[Theorem]{Corollary}
\newtheorem{Lemma}[Theorem]{Lemma}
\newtheorem{Proposition}[Theorem]{Proposition}
 { \theoremstyle{definition}
\newtheorem{Definition}[Theorem]{Definition}

\newtheorem{Example}[Theorem]{Example}
\newtheorem{Remark}[Theorem]{Remark}
\newtheorem{notn}[Theorem]{Notation}}

\newtheorem{case}{Case}

\newtheorem{subcase}{Case}

\newcommand{\la}{\lambda}
\newcommand{\al}{\alpha}
\newcommand{\be}{\beta}

\newcommand{\om}{\omega}

\newcommand{\ep}{\epsilon}
\newcommand{\De}{\Delta}

\newcommand{\si}{\sigma}

\newcommand{\R}{{\mathbb R}}
\newcommand{\Z}{{\mathbb Z}}

\newcommand{\abs}[1]{\left|#1\right|}
\newcommand{\deriv}[2]{\frac{\partial #1}{\partial #2}}
\newcommand{\C}{{\mathbb C}}
\renewcommand{\geq}{\geqslant}
\renewcommand{\leq}{\leqslant}

\newcommand{\varep}{\varepsilon}

\newcommand{\vect}[2]{\left(\begin{matrix}#1\\#2\end{matrix}\right)}
\newcommand{\matr}[4]{\left(\begin{matrix}#1&#2\\#3&#4\end{matrix}\right)}
\newcommand{\sltz}{{{\rm SL}_2(\Z)}}
\newcommand{\psl}{{{\rm PSL}_2(\Z)}}
\newcommand{\sltr}{{{\rm SL}_2(\R)}}

\begin{document}

\allowdisplaybreaks

\newcommand{\arXivNumber}{1502.07698}

\renewcommand{\PaperNumber}{016}

\FirstPageHeading

\ShortArticleName{Classifying Toric and Semitoric Fans by Lifting Equations from ${\rm SL}_2({\mathbb Z})$}

\ArticleName{Classifying Toric and Semitoric Fans\\ by Lifting Equations from $\boldsymbol{{\rm SL}_2({\mathbb Z})}$}

\Author{Daniel M.~KANE, Joseph PALMER and \'Alvaro PELAYO}
\AuthorNameForHeading{D.M.~Kane, J.~Palmer and \'A.~Pelayo}
\Address{University of California, San Diego, Department of Mathematics, \\
 9500 Gilman Drive \#0112, La Jolla, CA 92093-0112, USA}
\Email{\href{mailto:dakane@ucsd.edu}{dakane@ucsd.edu}, \href{mailto:j.palmer@math.rutgers.edu}{j.palmer@math.rutgers.edu}, \href{mailto:alpelayo@ucsd.edu}{alpelayo@ucsd.edu}}

\ArticleDates{Received April 17, 2017, in f\/inal form February 13, 2018; Published online February 22, 2018}

\Abstract{We present an algebraic method to study four-dimensional toric varieties by lifting matrix equations from the special linear group ${\rm SL}_2({\mathbb Z})$ to its preimage in the universal cover of ${\rm SL}_2({\mathbb R})$. With this method we recover the classif\/ication of two-dimensional toric fans, and obtain a description of their semitoric analogue. As an application to symplectic geometry of Hamiltonian systems, we give a concise proof of the connectivity of the moduli space of toric integrable systems in dimension four, recovering a known result, and extend it to the case of semitoric integrable systems with a f\/ixed number of focus-focus points and which are in the same twisting index class. In particular, we show that any semitoric system with precisely one focus-focus singular point can be continuously deformed into a system in the same isomorphism class as the Jaynes--Cummings model from optics.}


\Classification{52B20; 15B36; 53D05}

\rightline{{\it To Tudor S.~Ratiu on his $\text{65}^{\text{th}}$ birthday, with admiration.}}

\section{Introduction}

Toric varieties~\cite{Cox2003, CoLiSc2011, Da1978, Fultontoric, MiToric2008, Oda1988} have been extensively studied in algebraic and dif\/ferential geometry and so have their symplectic analogues, usually called symplectic toric manifolds or toric integrable systems~\cite{daSi2008}. The relationship between symplectic toric manifolds and toric varieties has been understood since the 1980s, see for instance Delzant~\cite{De1988},
Guillemin~\cite{GuKaehler1994, GuMoment1994}. The article~\cite{DuPe2009} contains a coordinate description of this relation.

In this paper we present an algebraic viewpoint to study non-singular complete four-di\-men\-sio\-nal toric varieties, based on the study of matrix relations in the special linear group $\sltz$. Indeed, one can associate to a rational convex polygon $\De$ the collection of primitive integer inwards pointing normal vectors to its faces, called a \emph{toric fan}. This is a $d$-tuple
\begin{gather*}
 (v_0=v_d, v_1, \ldots, v_{d-1})\in\big(\Z^2\big)^{d},
\end{gather*}
where $d\in\Z$ is the number of faces. A \emph{Delzant polygon}
(or \emph{toric polygon}) is one for which $\det(v_i, v_{i+1}) = 1$
for each $0\leq i \leq d-1$. This determinant condition forces the
vectors to satisfy the linear equations
$a_i v_i = v_{i-1}+v_{i+1}$, for $i=0, \ldots, d-1$ where $v_{-1}=v_{d-1}$,
which are parameterized
by integers $a_0, \ldots, a_{d-1}\in\Z$ (see Lemma~\ref{lem_av} and Fulton~\cite[p.~43]{Fultontoric}). These integers satisfy
\begin{gather}\label{eqn_intromatr}
 \matr{0}{-1}{1}{a_0} \matr{0}{-1}{1}{a_{1}}\cdots\matr{0}{-1}{1}{a_{d-1}} = \matr{1}{0}{0}{1}
\end{gather}
(this equation appears in \cite[p.~44]{Fultontoric}) but in fact, as can be seen in the exercises on p.~44 of~\cite{Fultontoric},
not all integers satisfying equation \eqref{eqn_intromatr} correspond to a toric fan. To deal with this issue we take an approach similar to that employed in~\cite{PR2000}. In particular, the identity in equation~\eqref{eqn_intromatr} in $\sltz$ can be lifted to the group $G$ presented as
\begin{gather*}
 G\cong \big\langle S,T \,|\, STS = T^{-1}ST^{-1}\big\rangle
\end{gather*}
and which, we will see in Lemma~\ref{lem_sltzpresn}, satisf\/ies $\sltz \cong G / \big(S^4\big)$. The group $G$, as is shown in \cite[Section~8.3]{PR2000}, is isomorphic to the pre-image of $\sltz$ in the universal cover of $\sltr$ (Proposition~\ref{prop_universalcover}), and thus we can def\/ine what we call the \emph{winding number} of an element of $g\in G$ that evaluates to the identity in $\sltz$. Roughly speaking, we view $g$ as a word in $S$ and $T$ and by applying this word to a vector one term at a time we produce a path around the origin. We def\/ine the winding number of $g$ to be the winding number of this path in the classical sense.

By considering a lift of equation~\eqref{eqn_intromatr} to the group~$G$, we obtain an equation with the property that a collection of integers $a_0, \ldots, a_{d-1}$ is a solution if and only if they correspond to a toric fan. Furthermore, from the integers it is straightforward to recover
the fan up to the appropriate isomorphism. Thus, the collection of toric fans can be studied by instead studying all $a_0, \ldots, a_{d-1}$ which
satisfy the equation in $G$. This allows us to prove some classical results about toric fans from a new perspective (Section~\ref{sec_toric}), and generalize these results to the semitoric case (Section~\ref{sec_semitoric}).

Recall from \cite{PeVNsemitoricinvt2009, PeVNconstruct2011} that a semitoric integrable system is a completely integrable system with two degrees of freedom and for which one of the Hamiltonians generates a periodic f\/low (see Def\/inition~\ref{def_stsystem}). Semitoric systems are a~generalization of toric integrable systems which can have, in addition to elliptic type singularities, also focus-focus singularities
(called nodal singularities in the context of Lefschetz f\/ibrations and algebraic geometry). Focus-focus singularities appear in algebraic
geometry~\cite{GrSi2006}, symplectic topology~\cite{ElPo2010,LeSy2010, Sy2003, Vi2014}, and many simple physical models such as the spherical pendulum~\cite{AbMa1978} and the Jaynes--Cummings system~\cite{Cu1965, JaCu1963, PeVNspin2012}. Associated to a~semitoric system there is also a collection of vectors $(v_0, \ldots, v_{d-1})$, def\/ined up to the appropriate notion of isomorphism, which satisfy more complicated equations (given explicitly in Def\/inition \ref{def_stfan}). Any collection of such vectors is henceforth referred to as a \emph{semitoric fan}. A semitoric fan can be thought of as a toric fan for which the relations between some pairs of adjacent vectors have been changed as a result of the presence of the focus-focus singularities. Roughly speaking, semitoric fans encode aspects of the singular af\/f\/ine structure induced by the singular f\/ibration associated to a semitoric integrable system. This af\/f\/ine structure also plays a~role in parts of symplectic topology, see for
instance Borman--Li--Wu~\cite{BoLiWu2014}, and mirror symmetry, see Kontsevich--Soibelman~\cite{GrSi2010,GrSi2011,KoSo2006}.

We present two theorems in this paper which are applications of the algebraic method we introduce. The f\/irst gives a classif\/ication of semitoric fans and the second is an application of the f\/irst which describes the path-connected components of the moduli space of semitoric integrable systems with a f\/ixed number of focus-focus singular points.

\begin{Theorem}\label{thm_A} Any semitoric fan may be obtained from a standard semitoric fan in a finite number of steps using four standard
 transformations.
\end{Theorem}

A more detailed description of Theorem \ref{thm_A} is given in Theorem~\ref{thm_stpoly} and the def\/initions of the standard semitoric fans and the four standard transformations used in Theorem~\ref{thm_A} are given in Def\/inition~\ref{def_fantrans}. As an application of Theorem~\ref{thm_A} to the symplectic geometry of integrable systems we will prove Theorem~\ref{thm_B}, which is the second main theorem of the paper.

The articles \cite{PeVNsemitoricinvt2009, PeVNconstruct2011} give a classif\/ication of semitoric systems in terms of several invariants, which can be essentially encoded into a polygon (a generalization of its toric counterpart and very closely related to a semitoric fan) together
with f\/initely many marked interior points, each of which is labeled by a formal power series in two variables and an integer. In contrast to the toric case, a semitoric system has a~family of polygons associated to it, each constructed using a toric momentum map def\/ined by removing certain closed subsets; the collection of all such polygons associated to a given system is henceforth referred to as the \emph{semitoric
polygon} associated to the system. The marked interior points are in fact the images of the focus-focus points of the integrable system
and the formal power series in two variables determines the semiglobal model around the focus-focus f\/iber up to a suitable notion of isomorphism (cf.~\cite{VNSnotes,VN2003}). The collection of integers labeling the focus-focus points for each polygon is the \emph{twisting index invariant}. This invariant encodes how the semiglobal model of the focus-focus f\/iber ``sits'' relative to the toric momentum map used to def\/ine each element of the semitoric polygon, so the integer assigned to a specif\/ic focus-focus point depends on the choice of polygon but the dif\/ference between the twisting index of consecutive (ordered according to the $x$-component) values does not depend on the choice of polygon and, roughly speaking, this dif\/ference measures the ``twist'' in the topology of the singular Lagrangian f\/ibration between them. Here we have not described the group action on the space of possible invariants, which is necessary so that choices made when constructing the invariants do not ef\/fect the outcome. See Section~\ref{sec_invariantsofst} for the precise def\/inition of these invariants and for a discussion of the group action.

For the next result we need to consider the moduli space of semitoric systems as a topological space so we can study continuous families in this space. The natural topology in this situation is the one induced by the metric on the moduli space of semitoric systems~\cite{PaSTMetric2015} and
which is described in Section~\ref{sec_fullmetric}.\footnote{For this paper we only def\/ine this metric on certain subsets of the full moduli space because this is all that is necessary to produce the same topology, see Remark~\ref{rmk_metricff}.} The metric is produced by pulling back a metric
from the space of all invariants of semitoric integrable systems. The role played by the polygon invariant in the def\/inition of this metric is related to the Duistermaat--Heckman measure.

The twisting index invariant of a system is an equivalence class of a list of integers assigned to each choice of toric momentum map. The integer associated to a given focus-focus point is not a well-def\/ined invariant when considered individually, but one can def\/ine the ``twisting index class" of a system. Two systems are in the same twisting index class if there exists a choice of polygon for each system such that the integers for these polygons agree (the precise notion is in Def\/inition~\ref{def_twistingindexclass}). Now we can state our second main result.

\begin{Theorem}\label{thm_B} With respect to the topology induced by the metric defined in~{\rm \cite{PaSTMetric2015}}, any two semitoric systems with the same number of focus-focus singular points which are in the same twisting index class may be continuously deformed into one another, up to isomorphism, via a continuous path of semitoric systems with the same number of focus-focus singular points and in the same twisting index class.
\end{Theorem}

We state a detailed version of this result later in the paper as Theorem \ref{thm_stconnect}.

\begin{Remark}The topology from~\cite{PaSTMetric2015} automatically places isomorphism classes of systems with dif\/ferent numbers of focus-focus points into dif\/ferent components of the space. It would be very interesting to have a~version of the metric, or at least topology, from \cite{PaSTMetric2015} which would allow us to study deformations of an isomorphism class of semitoric systems into another with a dif\/ferent number of focus-focus points but such a comparison gives rise to several nontrivial issues which are beyond the scope of the present paper. The most pressing issue being if such a deformation can really occur through semitoric systems, i.e., whether the non-degeneracy condition can be satisf\/ied for all values of the deformation (under some extra natural hypotheses, degeneracies must occur, see~\cite[Proposition~2.8]{HP2018}).
\end{Remark}

The Jaynes--Cummings system is an important example of a semitoric system with precisely one focus-focus point and is studied for example
in \cite{PeVNspin2012} (systems with exactly one focus-focus point are referred to
as systems of \emph{Jaynes--Cummings type} for this reason).
Since there is only one twisting index class of semitoric systems with exactly
one focus-focus point, Theorem~\ref{thm_B} implies the following.

\begin{Corollary} With respect to the topology induced by the metric defined in~{\rm \cite{PaSTMetric2015}}, any isomorphism class of semitoric system with precisely one focus-focus singular point may be continuously deformed into the isomorphism class of the Jaynes--Cummings system via
 a continuous path of isomorphism classes of semitoric systems, each of which have precisely one focus-focus point.
\end{Corollary}

Semitoric systems have been studied by mathematicians and physicists in the past decade, and there have been contributions to their study from many angles, including mathematical physics (see, e.g., Babelon--Dou{\c{c}}ot \cite{BaDo2012,BaDo2015}, Dullin \cite{Du2013}). While this paper deals with classical integrable systems, much of the work on these systems is motivated by inverse spectral problems about quantum integrable systems as pioneered in the work of Colin de Verdi{\`e}re \cite{CV1979,CV1980} and others, and which also has been the subject of recent works~\cite{ChPeVN2013,Ze1998}.

The paper is divided into two blocks. The f\/irst one concerns toric and semitoric fans and requires no prior knowledge of symplectic or algebraic geometry, while the second block, which consists only of Section~\ref{sec_stsyst}, contains applications to symplectic geometry and will probably be most interesting to those working on dif\/ferential geometry or Hamiltonian systems. The structure of the paper is as follows. In Section~\ref{sec_results} we state our main results, and the applications to symplectic geometry. In Section~\ref{sec_algsetup} we def\/ine the
necessary algebraic structures and prove several general algebraic results. In Section~\ref{sec_toric} we use this new algebraic approach to recover classical results about toric fans and in Section~\ref{sec_semitoric} we generalize these results to semitoric fans. Indeed, in Section~\ref{sec_stsyst} we use the results of Section~\ref{sec_semitoric} to study the connectivity of the moduli space of semitoric systems.

\section{Fans, symplectic geometry, and winding numbers}\label{sec_results}

\subsection{Toric fans}
A toric variety is a variety which contains an algebraic torus as a~dense (in the Zariski topology) open subset such that the standard action of the torus on itself can be extended to the whole variety. That is, a toric variety is the closure of an algebraic torus orbit~\cite{Fultontoric}. By an
algebraic torus we mean the product $\C^*\times\cdots\times\C^*$, where $\C^* = \C \setminus \{0\}$. Under some mild assumptions (which are automatically satisf\/ied if the variety is smooth), the geometry of a~toric variety is completely determined by the associated fan~\cite[Section 1.5]{Fultontoric}.
\begin{Definition}\label{def_generalfan}
A \emph{rational strongly convex cone} is a convex cone of a vector space with apex at the origin generated by a f\/inite number of integral vectors which contains no line through the origin. A \emph{fan} is a set of rational strongly convex cones in a real vector space such that the face of each cone is also a cone and the intersection of any two cones is a face of each.
\end{Definition}

Def\/inition~\ref{def_generalfan} is the general notion of a fan in arbitrary dimension; in this paper we will be concerned with two-dimensional smooth compact toric varieties and their associated fans, which for simplicity we will call toric fans. We identify such fans with a sequence of points in~$\Z^2$.
\begin{Definition}
A \emph{toric fan} is a sequence of lattice points \begin{gather*}(v_0=v_d, v_1, \ldots, v_{d-1})\in \big(\Z^2\big)^d\end{gather*} labeled in
counter-clockwise order such that each pair of adjacent vectors generates all of $\Z^2$ and the angle between any two adjacent vectors is less than $\pi$ radians. That is, $\det(v_i, v_{i+1}) = 1$ for $i = 0, \ldots, d-1$.
\end{Definition}

\begin{Definition}\label{def_poly} For the purposes of this paper by \emph{convex polygon} we mean the intersection in~$\R^2$ of f\/initely or inf\/initely many closed half planes such that there are at most f\/initely many corner points in each compact subset of $\R^2$.
\end{Definition}

\begin{Definition}\label{def_delzantpoly}
A \emph{Delzant polygon} (or \emph{toric polygon}) is a compact convex polygon $\De$ in $\R^2$ which is simple, rational, and smooth. Recall $v\in\Z^2$ is primitive if $v=kw$ for some $k\in\Z_{>0}$ and $w\in\Z^2$ implies $k=1$ and $w=v$.
\begin{enumerate}\itemsep=0pt
 \item[1)] $\De$ is \emph{simple} if there are exactly two edges meeting at each vertex;
 \item[2)] $\De$ is \emph{rational} if for each edge $e$ of $\De$ there exists a vector in $\Z^2$ which is normal to $e$;
 \item[3)] $\De$ is \emph{smooth} if the inwards pointing primitive vectors normal to any pair of adjacent edges form a basis of $\Z^2$.
\end{enumerate}
\end{Definition}

Delzant polygons were introduced in the work of Delzant~\cite{De1988} in symplectic geometry, who built on the work of Atiyah~\cite{At1982}, Kostant~\cite{Ko1973}, and Guillemin-Sternberg~\cite{GuSt1982} to give a classif\/ication of symplectic toric $4$-manifolds in terms of the Delzant polygon (in fact, their work was in any dimension and higher dimensional toric integrable systems correspond to higher dimensional \emph{Delzant polytopes} which satisfy analogous conditions to those given in this Def\/inition~\ref{def_delzantpoly}). Delzant polygons are similar to Newton polygons, as in~\cite{Fultontoric}, except that the vertices of Delzant polygons do not have to have integer coordinates. Just as in the case with Newton polygons, a toric fan may be produced from a Delzant polygon by considering the collection of inwards pointing normal vectors of the polygon. Notice that since Delzant polygons are required to be convex we automatically have that the angle between adjacent vectors in the fan is less than~$\pi$. The natural action of $\sltz$ on $\R^2$ induces an action on the set of Delzant polygons. In \cite{Fultontoric} Fulton stated the following result for toric fans. We have adapted the statement to relate it to Delzant polygons.

\begin{Theorem}[{Fulton~\cite[p.~44]{Fultontoric}}]\label{thm_fultonpolygon} Up to the action of $\sltz$, every Delzant polygon can be obtained from a Delzant triangle, rectangle, or Hirzebruch trapezoid by a finite number of corner chops.
 \end{Theorem}

 \begin{figure}[t]\centering
 \includegraphics[width=340pt]{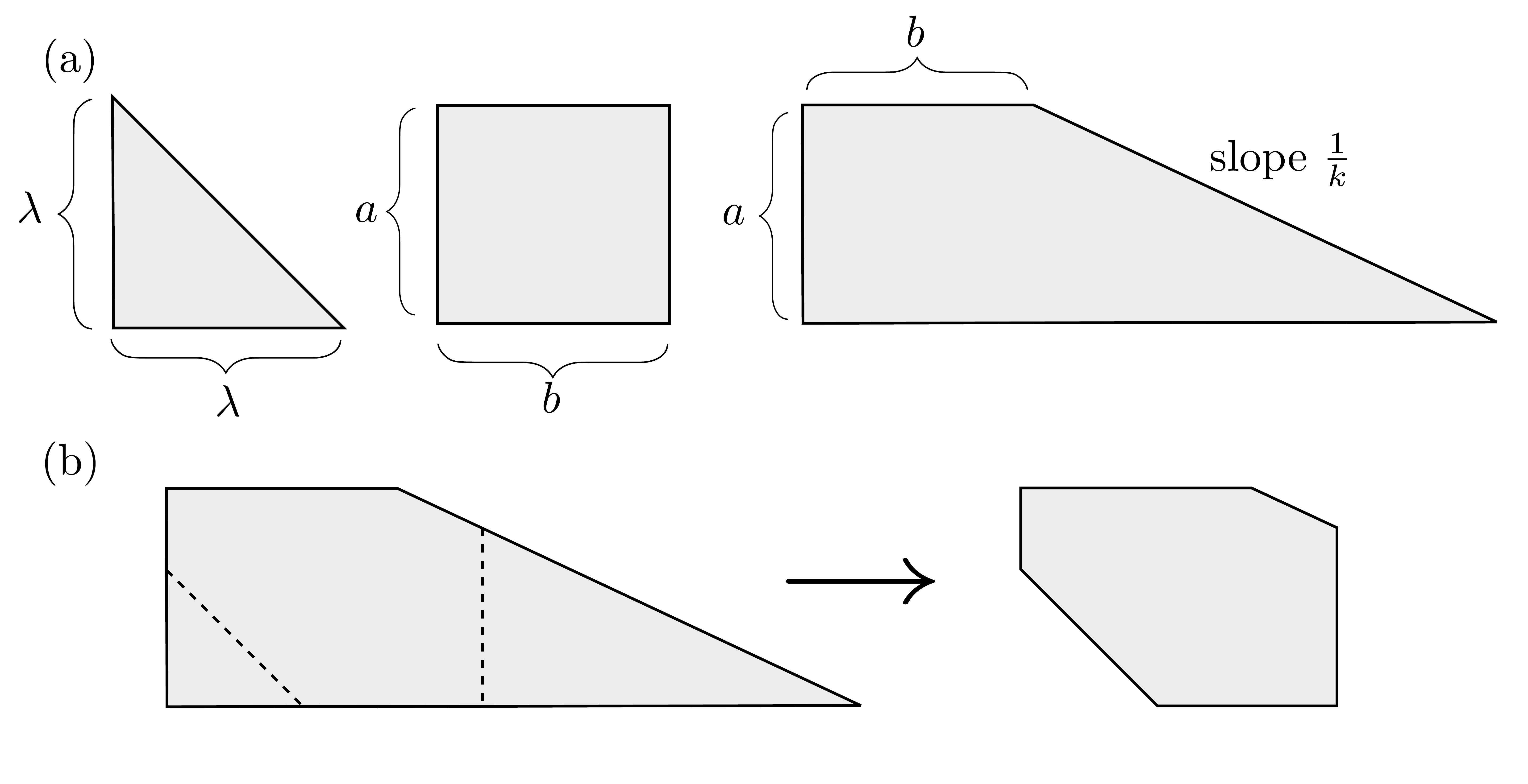}
 \caption{(a) The three minimal models from Theorem \ref{thm_fultonpolygon}. (b)~An illustration of a Delzant polygon produced by corner chopping the Hirzebruch trapezoid.} \label{fig_minimalmodels}
 \end{figure}

The minimal models (the Delzant triangle, rectangle, and Hirzebruch trapezoid with parameter $k>1$) are def\/ined in Def\/inition~\ref{def_trisquaretrap} and depicted in Fig.~\ref{fig_minimalmodels}. The corner chop operation is def\/ined in Def\/inition~\ref{def_fantrans} and is a standard operation in algebraic and symplectic geometry (see, for instance,~\cite{KaKe2007}). It corresponds to an equivariant symplectic blowup.

The proof of Theorem \ref{thm_fultonpolygon} sketched by Fulton in~\cite[Section~2.5]{Fultontoric}, which uses only two-dimensional geometry and basic combinatorial arguments, is relatively long and does not immediately generalize to the case of semitoric polygons we study in this paper. In~Section \ref{sec_toric} we provide an alternative proof using $\sltz$-relations which may be easily extended to the semitoric case.

The moduli space of toric polygons was endowed with a metric given by the Lebesgue measure of the symmetric dif\/ference in~\cite{PePRS2013}, which we denote by $d_\mathbb{T}$, and in that paper the following consequence of Theorem~\ref{thm_fultonpolygon} was proved:

\begin{Theorem}[\cite{PePRS2013}]\label{thm_toricconnected}
 The moduli space of toric polygons is path-connected with respect to the topology induced by $d_\mathbb{T}$.
\end{Theorem}

That is, any two toric polygons may be deformed onto each other continuously via a path of toric polygons. One shows this by f\/irst knowing how to generate all toric polygons as in Theorem~\ref{thm_fultonpolygon}. Then one shows, using elementary analysis, that the three minimal models can be continuously transformed into one another and that the corner chop operation is continuous. Again, let us emphasize that Theorem~\ref{thm_fultonpolygon} is not new, but in this paper we will prove it from a~new viewpoint which we believe to be more natural and which generalizes to the case of semitoric fans.

\subsection{Semitoric fans}

Motivated by semitoric polygons (originally def\/ined in \cite[Def\/inition 2.5]{PeVNconstruct2011}) we def\/ine semitoric fans. A~semitoric polygon is a family of polygons which are similar to Delzant polygons except that the presence of the focus-focus singular points in semitoric systems may cause some corners of the polygons to not satisfy the smoothness condition of Def\/inition~\ref{def_delzantpoly}. Such corners must satisfy other conditions and are known as fake or hidden corners. They can exist on the top or bottom boundary of the polygon but there will always be an inf\/inite subfamily of semitoric polygons related to the system for which the hidden and fake corners are all on the top boundary, and we will use one of these polygons to produce the associated semitoric fan. The choice of which polygon with this property is used is ref\/lected in the symmetry group of semitoric fans which we introduce in Def\/inition~\ref{def_stsymmetry}. Let
\begin{gather*}
 T = \matr{1}{1}{0}{1}.
\end{gather*}

\begin{Definition}[{\cite[Section 4.1 and Def\/inition 4.1]{PeVNconstruct2011}}]\label{def_cornerconditions} Let $v,w\in\Z^2$. The ordered pair $(v,w)$ of vectors:
 \begin{enumerate}\itemsep=0pt
 \item[1)] is \emph{on the top boundary} if both vectors are in the open lower half-plane;
 \item[2)] satisf\/ies the \emph{Delzant condition} if $\det(v, w)=1$;
 \item[3)] satisf\/ies the \emph{hidden condition} if $\det(v, Tw)=1$; and
 \item[4)] satisf\/ies the \emph{fake condition} if $\det(v, Tw)=0$.
 \end{enumerate}
\end{Definition}

\begin{Remark}\label{rmk_fakeanddelzant} Notice that a pair $(v,w)$ satisf\/ies both the fake and Delzant conditions if and only if
\begin{gather*}
 v = \vect{k+\varepsilon}{\varepsilon} \qquad \textrm{and} \qquad w = \vect{k}{\varepsilon}
\end{gather*}
for some $k\in\Z$ and $\varepsilon\in\{-1,+1\}$ and in order for such a~pair to be in the top boundary we can only have the case in which
 $\varepsilon = -1$.
\end{Remark}

\begin{Definition}\label{def_stfan} Let $d\in\Z$ with $d>2$. A \emph{semitoric fan} is a collection of primitive vectors $(v_0=v_d, v_1,\ldots,v_{d-1})\in\big(\Z^2\big)^d$ labeled in counter-clockwise order such that each pair of adjacent vectors $(v_i, v_{i+1})$ for $i\in\{0,\ldots,d-1\}$ is labeled as a \emph{Delzant}, \emph{fake}, or \emph{hidden} corner. We require that each labeled pair of vectors satisf\/ies the corresponding condition from Def\/inition~\ref{def_cornerconditions} and we further require that all fake and hidden corners be on the top boundary. The \emph{defect} of a~semitoric fan is the number of corners which are either
 fake or hidden.
\end{Definition}

Notice that the labeling of the pairs is required only because of the case described in Remark~\ref{rmk_fakeanddelzant} in which a pair can satisfy both the fake and Delzant conditions. In all other cases the corner type of a pair of vectors can be uniquely determined by inspecting the
vectors involved.

Def\/inition \ref{def_stfan} is inspired by the toric case. Theorem \ref{thm_fultonpolygon} states that any toric fan can be produced from a minimal model using only corner chops. Similarly, our goal is to use a series of transformations to relate any semitoric fan to a standard form up to the action of the appropriate symmetry group.

\begin{Definition}\label{def_stsymmetry}
 The \emph{symmetry group of semitoric fans} is given by
 \begin{gather*}
 \mathcal{G}' = \big\{ T^k \,|\, k\in\Z\big\},
 \end{gather*}
 where $\mathcal{G}'$ acts on a semitoric fan by acting on each vector in the fan in the standard fashion.
\end{Definition}

\begin{Definition}\label{def_fantrans}\quad
 \begin{enumerate}\itemsep=0pt
 \item Let $c\in\Z_{\geq0}$. The \emph{standard semitoric fan of defect $c$} is the fan $(u_0, \ldots, u_{c+3})\in\big(\Z^2\big)^{c+4}$ given by
 \begin{gather*}
 u_0 = \vect{0}{-1}, \qquad u_1 = \vect{1}{0},\qquad u_2 = \vect{c}{1}, \qquad u_3 = \vect{-1}{0},
 \end{gather*}
 and
 \begin{gather*}
 u_{4+n} = \vect{-c+n}{-1}
 \end{gather*}
 for $n = 0, \ldots, c-1$ in which the f\/irst four pairs of vectors are Delzant corners and the rest are fake corners.
 \item Let $(v_0=v_d, \ldots, v_{d-1})\in\big(\Z^2\big)^d$ be a semitoric fan. The following are called the \emph{four fan transformations}:
 \begin{enumerate}\itemsep=0pt
 \item\label{trans_chop} Suppose that $(v_i, v_{i+1})$ is a Delzant corner for some $i\in\{0, \ldots, d-1\}$. Then
 \begin{gather*}
 (v_0, \ldots, v_i, v_i + v_{i+1}, v_{i+1}, \ldots, v_{d-1})\in\big(\Z^2\big)^{d+1}
 \end{gather*}
 obtained by inserting the sum of two adjacent vectors between them is a semitoric fan with the new pairs $(v_i, v_i+v_{i+1})$ and $(v_i+v_{i+1}, v_{i+1})$ both labeled as Delzant corners (Lemma~\ref{lem_cc}). The process of producing this new fan from the original generalizes \emph{corner chopping} for toric fans (cf.~\cite{KaKe2007}). For this reason, it is also referred to as a~corner chop(ping) in the context of semitoric fans.
 \item Suppose that $(v_{i-1},v_i)$ and $(v_i, v_{i+1})$ are Delzant corners and that $v_i = v_{i-1} + v_{i+1}$. Then
 \begin{gather*}
 (v_0, \ldots, v_{i-1}, v_{i+1}, \ldots, v_{d-1})\in \big(\Z^2\big)^{d-1}
 \end{gather*}
 is a semitoric fan with the pair $(v_{i-1}, v_{i+1})$ being labeled as a Delzant corner (Lem\-ma~\ref{lem_reversecc}). This process of producing this new fan from the original generalizes \emph{reverse corner chopping} for toric fans (cf.~\cite{KaKe2007}). For this reason, it is also referred to as a~reverse corner chop(ping) in the context of semitoric fans.
 \item Suppose that the pair $(v_i, v_{i+1})$ is a~hidden corner. Then
 \begin{gather*}
 (v_0, \ldots, v_i, Tv_{i+1}, v_{i+1}, \ldots, v_{d-1})\in\big(\Z^2\big)^{d+1}
 \end{gather*}
 is a semitoric fan with $(v_i, Tv_{i+1})$ a Delzant corner and $(Tv_{i+1}, v_{i+1})$ a fake corner (Lemma~\ref{lem_removehidden}). The process of producing this fan is called \emph{removing the hidden corner $(v_i, v_{i+1})$}.
 \item Suppose that the pair $(v_i, v_{i+1})$ is a fake corner and the pair $(v_{i+1}, v_{i+2})$ is a~Delzant corner on the top boundary. Then
 \begin{gather*}
 (v_0, \ldots, v_i, Tv_{i+2}, v_{i+2}, \ldots, v_{d-1})\in\big(\Z^2\big)^{d}
 \end{gather*}
 is a semitoric fan with $(v_i, Tv_{i+2})$ a Delzant corner and $(Tv_{i+2}, v_{i+2})$ a fake corner (Lemma~\ref{lem_fakedelzcommute}). The process of producing this fan is called \emph{commuting a~fake and a~Delzant corner}.
 \end{enumerate}
 \end{enumerate}
\end{Definition}

\begin{Remark}The corner chop and reverse corner chop can be obtained as equivariant symplectic blowups and blowdowns, respectively.
\end{Remark}

\begin{Remark} The standard semitoric fans from Def\/inition~\ref{def_fantrans} are, for the purposes of this paper, analogous to the minimal models from toric geometry (as in Theorem~\ref{thm_fultonpolygon}) because any semitoric fan can be transformed into a standard semitoric fan using the four fan transformations. That being said, the set of standard semitoric fans does not correspond to the set of all \emph{minimal semitoric fans}, which are those fans which do not admit a reverse corner chop. The problem of classifying all minimal semitoric fans is being addressed in a future paper (which is now available, see~\cite{KaPaPe2018}).
\end{Remark}

The transformations of removing a hidden corner and commuting a fake and a Delzant corner can also be reversed, but we will not need the inverses of those operations for the proofs in this paper, so we do not list them here. Using the algebraic results from Section \ref{sec_algsetup} we show the following in Section~\ref{sec_semitoric}.

\begin{Theorem}\label{thm_stpoly}
Let $d\geq 3$ be an integer. Any semitoric fan $(v_0, \ldots, v_{d-1})\in\big(\Z^2\big)^d$ of defect $c\in\Z_{\geq0}$ may be transformed into a semitoric fan $\mathcal{G}'$-equivalent to the standard semitoric fan of defect $c$ by using the four fan transformations finitely many times.
\end{Theorem}

\begin{Remark} The method we are using to study semitoric manifolds is analogous to the method we use to study toric manifolds. Theorem~\ref{thm_fultonpolygon} explains how to generate the toric polygons and is used to prove that the space of toric polygons is path-connected (Theorem~\ref{thm_toricconnected}) which implies that the space of toric manifolds is path-connected (Theorem~\ref{thm_toricmanifoldconnect}). Similarly, Theorem~\ref{thm_stpoly} shows how to generate the semitoric polygons, and as an application we prove Lemma~\ref{lem_stpolyconnect} which describes the connected components in the space of semitoric ingredients (Def\/inition~\ref{def_listofingre}) and this implies Theorem~\ref{thm_stconnect}, which describes the connected components of the moduli space of semitoric systems.
\end{Remark}

\begin{Remark}It is likely that there are connections between our work and that of Gross and Siebert \cite{GrSi2006, GrSi2010, GrSi2011}.
In particular, it seems to us that one of the objects that Gross and Siebert use is very closely related to the notion of a semitoric fan. Understanding this would be an entirely dif\/ferent project which merits further study.
\end{Remark}

\subsection{Algebraic tools: the winding number}

It is shown in Lemma \ref{lem_sltzpresn} that the special linear group $\sltz$ may be presented as
\begin{gather*}
 \sltz = \big\langle S,T \,|\, T^{-1}ST^{-1} = STS,\ S^4 = I\big\rangle,
\end{gather*}
where
\begin{gather*}
 S = \matr{0}{-1}{1}{0}\qquad \text{and} \qquad T=\matr{1}{1}{0}{1}.
\end{gather*}
Thus equation~\eqref{eqn_intromatr} becomes
\begin{gather}\label{eqn_STtoric}
 ST^{a_0}\cdots ST^{a_{d-1}} = I,
\end{gather}
where $a_0, \ldots, a_{d-1}\in\Z$ and $I$ denotes the $2\times 2$ identity matrix. Given $v_0, v_1\in\Z^2$ with $\det(v_1, v_2)=1$ a~set of vectors \begin{gather*}(v_0, v_1, \ldots, v_{d-1})\in\big(\Z^2\big)^d\end{gather*} may be produced by \begin{gather*} v_{i+2} = -v_i + a_i v_{i+1}\end{gather*} for
$i=0, \ldots, d-1$ where we def\/ine $v_d = v_0$ and $v_{d+1}= v_1$. In this way associated to each list of integers satisfying equation~\eqref{eqn_intromatr} there is an ordered collection of vectors unique up to $\sltz$. It can be seen that the determinant between
any adjacent pair of these vectors is one and thus if these vectors are labeled in counter-clockwise order, then they are a toric fan.
The reason that not all sequences of integers which satisfy equation~\eqref{eqn_intromatr} correspond to a toric fan is that the vectors
$v_0, \ldots, v_{d-1}\in\Z^2$ may circle more than once around the origin, and thus not be labeled in counter-clockwise order (see Fig.~\ref{fig_windtwice}). Thus, we see that equation~\eqref{eqn_STtoric} is merely a necessary and not suf\/f\/icient condition for a sequence $(a_0, \ldots, a_{d-1})$ to correspond to a toric fan.

Let $K = \mathrm{ker}(\Z*\Z \to \sltz)$ where $\Z*\Z$ denotes the free group with generators\footnote{Starting in Section~\ref{sec_algsetup} we will adopt Notation~\ref{notationST} to dif\/ferentiate between the generators $S$ and $T$ in the dif\/ferent groups.}~$S$ and~$T$ and the map~$\Z*\Z\to\sltz$ is the natural projection. For any word in $K$ a sequence of vectors may be produced by letting the word act on a vector $v\in\Z^2$ one term at a time. We know that this sequence ends back at~$v$, but the sequence of vectors produced contains more information about the word. This sequence can be used to def\/ine a path in $\R^2\setminus\{(0,0)\}$ by considering the piecewise linear path between the ends of the
vectors, see Def\/inition~\ref{def_circleorigin}. Of particular interest, especially when studying toric and semitoric fans, is the \emph{winding number} of such a path. That is, the number of times the path, and hence the collection of vectors, circles the origin. This construction is explained in detail in Section~\ref{sec_toric}, and in particular Def\/inition~\ref{def_circleorigin} gives a precise def\/inition of the
number of times an ordered collection of vectors circles the origin. Let $w\colon \Z*\Z\to\Z$ be the unique homomorphism satisfying
\begin{gather*}
 w(S)= 3\qquad \text{and}\qquad w(T) = -1.
\end{gather*}
We f\/ind that if $\si\in K$ then $w(\si)$ is a multiple of 12 and $\nicefrac{w(\si)}{12}$ is the winding number associated to the word $\si$, see Lemma~\ref{lem_windingnumber}. We present the group
\begin{gather*}
 G= \big\langle S,T \,|\, STS = T^{-1}ST^{-1}\big\rangle
\end{gather*}
on which $w$ descends to a well-def\/ined homomorphism $w_G\colon G\to\Z$. In fact, $G$ is isomorphic to the pre-image of~$\sltz$ in the universal cover of~$\sltr$ (Proposition~\ref{prop_universalcover}). Thus, if $K'$ is the image of~$K$ projected to~$G$, then given some $g\in K'$ there
is an associated closed loop in~$\sltr$. The fundamental group of~$\sltr$ is~$\Z$ and the classical winding number of this loop in~$\sltr$ coincides with $\nicefrac{w_G(g)}{12}$. Finally, in Corollary~\ref{cor_toricfaneqn} we show that integers $a_0, \ldots, a_{d-1}\in\Z$ correspond to a toric fan if and only if the equality
\begin{gather*}
 ST^{a_0}\cdots ST^{a_{d-1}} = S^4
\end{gather*}
is satisf\/ied in $G$. This correspondence is the basis of our method to study toric and semitoric fans.

\begin{Remark}
We note that~\cite{PR2000} also needed to make use of this winding number for similar reasons. They def\/ine the homomorphism~$w$ ($\Phi$ in their notation) in Section~8.4, and come up with an interesting interpretation of it. While~$w(\gamma)$ can be interpreted as a winding number for~$\gamma$ that map to the identity in $\sltz$, for other values of $\gamma$, the interpretation is not so obvious.

Firstly, \cite{PR2000} come up with an interpretation for elements of the universal cover of $\sltr$. In particular, they write elements of the group as pairs $(M,[\gamma])$ where $M\in \sltr$ and~$[\gamma]$ a~homotopy class of paths in $\R^2-\{0\}$ from $b$ to $Mb$ for some f\/ixed basepoint $b$. The composition rule can then be written as $(M,[\gamma])(M',[\gamma']) = (MM',[\gamma^{M'}\gamma'])$ where $\gamma^{M'}$ denotes the path obtained by applying $M'$ pointwise to~$\gamma$, and the multiplication with $\gamma'$ denotes concatenation.

In this language, \cite{PR2000} f\/ind an interpretation of $w$ relating it to the modular form $\Delta$. In particular, if $M\in \sltz$ is given
by $M=\left(\begin{smallmatrix} a & b \\ c & d\end{smallmatrix}\right)$, then treating $M$ as a linear fractional transformation in the standard way, we have that
\begin{gather*}
\Delta(Mz) = (cz+d)^{12}\Delta(z).
\end{gather*}
Taking logs of this one f\/inds that $\log(\Delta(Mz)) - \log(\Delta(z))$ is some branch of $12\log(cz+d)$. If~$M$ is given not just as a matrix, but as the endpoint of a path from the identity in $\sltr$, then this def\/ines the appropriate branch. If $M$ corresponds to the element $(M,[\gamma])$ in the universal cover of $\sltr$ (where $\gamma$ uses~$z$ as its basepoint), \cite{PR2000} show that
\begin{gather}\label{eqn_PhiDefine}
\log(\Delta(M z)) - \log(\Delta(z)) = 12 \log(cz+d) + 2\pi i w((M,\gamma)),
\end{gather}
where the branch of $\log(cz+d)$ is def\/ined using the path $\gamma$. In fact, \cite{PR2000} use equation~\eqref{eqn_PhiDefine} to def\/ine~$w$, and then show in Section~8.5 the action on generators. We note that from this expression, it is clear that if $M=I$, that since $\log(\Delta)$ is analytic on the upper half plane, that $w((M,\gamma))$ comes entirely from the branch of $\gamma z$. In particular, this means that $w((I,\gamma))$ is $12$ times the winding number of~$\gamma z$.
\end{Remark}

\subsection{Applications to symplectic geometry}

While Delzant polygons are in correspondence with closed toric manifolds, semitoric polygons are associated with the so-called semitoric integrable systems. A \emph{semitoric integrable system} (or \emph{semitoric manifold}) is given by a triple $(M, \om, F\colon M\to\R^2)$ where $(M,\om)$ is
a connected symplectic $4$-manifold and $F\colon M\to\R^2$ is an integrable system given by two maps $J,H\colon M\to\R^2$ such that $J$ is a proper map which generates a~periodic f\/low (see Def\/inition~\ref{def_stsystem} for the precise def\/inition). In~\cite{PeVNconstruct2011} the authors prove a result analogous to Delzant's, classifying semitoric systems satisfying a mild assumption via a list of ingredients which includes a family of polygons (there is an overview of this result in Section~\ref{sec_invariantsofst}).

In \cite{PaSTMetric2015} the second author def\/ines a metric space structure for the moduli space of semitoric systems using~\cite{PeVNconstruct2011}, which is related to the Duistermaat--Heckman measure. The metric is designed to induce a topology that respects continuous transformations of the invariants of semitoric systems (this follows tautologically from the def\/inition) and
to be compatible with the known topology on toric integrable systems from~\cite{PePRS2013}, in the sense that the projection map from isomorphism classes of compact semitoric systems with zero focus-focus points to toric systems is continuous, as proven in~\cite[Corollary~4.14]{PaSTMetric2015} (this map is an projection since the isomorphism for toric integrable systems is strictly weaker than that for semitoric systems). A~natural question is whether, with respect to a such structure, the space is path-connected. That is, can two isomorphism classes of semitoric systems be continuously
deformed (with respect to the above topology) into one another, via a continuous path of isomorphism classes of semitoric systems?

This is preceded by~\cite{PePRS2013} in which the authors construct a natural metric on the moduli space of symplectic toric manifolds which is related to the Duistermaat--Heckman measure and prove Theorem~\ref{thm_toricconnected}, which is used to conclude the following.

\begin{Theorem}[\cite{PePRS2013}]\label{thm_toricmanifoldconnect} With respect to the topology defined in~{\rm \cite{PePRS2013}}, the moduli space of toric manifolds is path-connected.
\end{Theorem}

Similarly, Theorem \ref{thm_stpoly} implies the following statement.

\begin{Theorem}\label{thm_stconnect}
 If $(M,\om,F)$ and $(M',\om',F')$ are simple semitoric integrable systems such that:
 \begin{enumerate}[label=(\roman*),font = \normalfont]\itemsep=0pt
 \item they have the same number of focus-focus singularities;
 \item they are in the same twisting index class,
 \end{enumerate}
 then there exists a continuous $($with respect to the topology defined in {\rm \cite{PaSTMetric2015})} path of semitoric systems with the same number of focus-focus points and in the same twisting index class between them, up to isomorphism. That is, the space of isomorphism classes of semitoric systems with fixed number of focus-focus points and in a fixed twisting index class is path-connected.
\end{Theorem}

The number of focus-focus singularities is discussed in Section~\ref{sec_numberffpoints}, the twisting index class of a~semitoric system is def\/ined in Def\/inition~\ref{def_twistingindexclass}, and the notion of a simple semitoric integrable system is def\/ined in Def\/inition~\ref{def_stsystem}.
\begin{Remark} It is important to note that while there is only one twisting index class for systems with exactly one focus-focus point, this does
not mean that the twisting index invariant is trivial in this case. That is, there exist systems with exactly one focus-focus point for which all invariants agree except for the twisting index and the systems are not isomorphic.
\end{Remark}

\section[Algebraic set-up: matrices and $\sltz$ relations]{Algebraic set-up: matrices and $\boldsymbol{\sltz}$ relations}\label{sec_algsetup}

\begin{notn}\label{notationST} We will present several dif\/ferent groups on generators $S$ and $T$ because it is important to be able to easily see the standard homomorphisms between these groups. When referring to a word in $S$ and $T$ we will use a~subscript to indicate which group that word belongs to. For instance, we write $(S^4)_\sltz$ to refer to the element of~$\sltz$. To denote the dif\/ferent equalities in these groups we will use an equal sign with the group in question as a subscript. That is, if an equality holds in the group $H$ we will write $=_H$, so $\si =_H \eta$ is shorthand for $(\si)_H = (\eta)_H$. For example, $S^4 =_{\sltz} I$ but $S^4 \neq_G I$. Finally, if $v\in\Z^2$ then $Tv$ will always mean $T_\sltz v$.
\end{notn}

The $2\times2$ special linear group over the integers, $\sltz$, is generated by the matrices
\begin{gather*}
 S_\sltz = \matr{0}{-1}{1}{0}\qquad \text{and}\qquad T_\sltz = \matr{1}{1}{0}{1}.
\end{gather*}
We will see that to each closed toric (resp.\ semitoric) integrable system there is an associated toric (resp.\ semitoric) fan and we will use the algebraic structure of $\sltz$ to study these fans. For our purposes, the following presentation of $\sltz$ will be the most natural way to view the group.

\begin{Lemma}\label{lem_sltzpresn}
The $2\times 2$ special linear group over the integers, $\sltz$, may be presented as
\begin{gather*}
 \sltz \cong \big\langle S,T \,|\, T^{-1}ST^{-1} = STS,\, S^4\big\rangle.
\end{gather*}
\end{Lemma}
\begin{proof}
 It is well-known that
 \begin{gather}\label{eqn_sltzwellknown}
 \sltz \cong \big\langle S,T \,|\, (ST)^3=S^2,\, S^4 \big\rangle
 \end{gather}
 (see for instance \cite[equation~(A.2)]{KaTu2008}).
 The result follows from equation~\eqref{eqn_sltzwellknown}
 and the observation that
 \begin{gather*}
 (ST)^3 =_\sltz S^2 \Leftrightarrow STS =_\sltz T^{-1}ST^{-1}.\tag*{\qed}
 \end{gather*}\renewcommand{\qed}{}
\end{proof}

For $v,w\in\Z^2$ let $[v,w]$ denote the $2\times2$ matrix with $v$ as the f\/irst column and $w$ as the second and let $\det(v,w)$ denote the determinant of the matrix $[v,w]$.

\begin{Lemma}\label{lem_av}
 Let $u,v,w\in\Z^2$ and $\det(u,v)=1$. Then $\det (v, w) = 1$ if and only if there exists some $a\in\Z$ such that $w = - u + av$.
\end{Lemma}
\begin{proof}
 In the basis $(u, v)$ we know that $v = \left(\begin{smallmatrix} 0\\ 1\end{smallmatrix}\right)$. Write $w = \left(\begin{smallmatrix} b\\ a\end{smallmatrix}\right)$ for some $a, b\in\Z$. Then we can see that $\det(v, w) = -b$ so $\det(v,w)=1$ if and only if $b = -1$. That is, $w = -u + a v$.
\end{proof}

The result of Lemma \ref{lem_av} can be easily summarized in a matrix equation, as we will now show. Let
\begin{gather*}
 (v_0 = v_d, v_1 = v_{d+1}, v_2,\ldots, v_{d-1}) \in\big(\Z^2\big)^d
\end{gather*}
be a toric fan and def\/ine $A_i = [v_i, v_{i+1}]$ for $i=0, \ldots, d $. Note that $A_d = A_0$.

\begin{Lemma}\label{lem_ast}
 For each $i\in 0, \ldots, d-1$ there exists an integer $a_i \in\Z$ such that $A_{i+1} = A_i S T^{a_i}$.
\end{Lemma}
\begin{proof}
 By the def\/inition of a toric fan we know that for each $0\leq i < d-2$ we have that
 \begin{gather*}
 \det(v_i, v_{i+1}) = \det(v_{i+1}, v_{i+2}) = 1
 \end{gather*}
 so by Lemma \ref{lem_av} there exists $a_i \in \Z$ such that $v_{i+2} = -v_i + a_i v_{i+1}$. Then
 \begin{gather*}
 A_i S T^{a_i} = [v_{i+1}, -v_i + a_i v_i]= [v_{i+1}, v_{i+2}]=A_{i+1},
 \end{gather*}
 and this concludes the proof.
\end{proof}

It follows that
\begin{gather*}
 A_d = A_{d-1}\big(ST^{a_{d-1}}\big)_{\sltz}= A_{d-2}\big(ST^{a_{d-2}}ST^{a_{d-1}}\big)_{\sltz} = \cdots\\
 \hphantom{A_d}{} = A_0 \big(ST^{a_0}\cdots ST^{a_{d-1}}\big)_\sltz,
\end{gather*}
which means $ A_0 = A_d = A_0 (ST^{a_0}\cdots ST^{a_{d-1}})_{\sltz}$,
and so
\begin{gather}\label{eqn_st}
ST^{a_0}\cdots ST^{a_{d-1}} =_{\sltz} I.
\end{gather}
This is a restatement of equation \eqref{eqn_intromatr} which is from \cite[p.~44]{Fultontoric}. So to each toric fan of $d$ vectors there is an associated $d$-tuple of integers which satisfy equation~\eqref{eqn_st}, but having a tuple of integers which satisfy equation~\eqref{eqn_st} is not enough to assure that they correspond to a~toric fan. The determinant of the vectors will be correct but, roughly speaking, if the vectors wind around the origin more than once then they will not be labeled in the correct order to be a~toric fan, as it occurs in the following example.

 \begin{figure}[t]\centering
 \includegraphics[width=200pt]{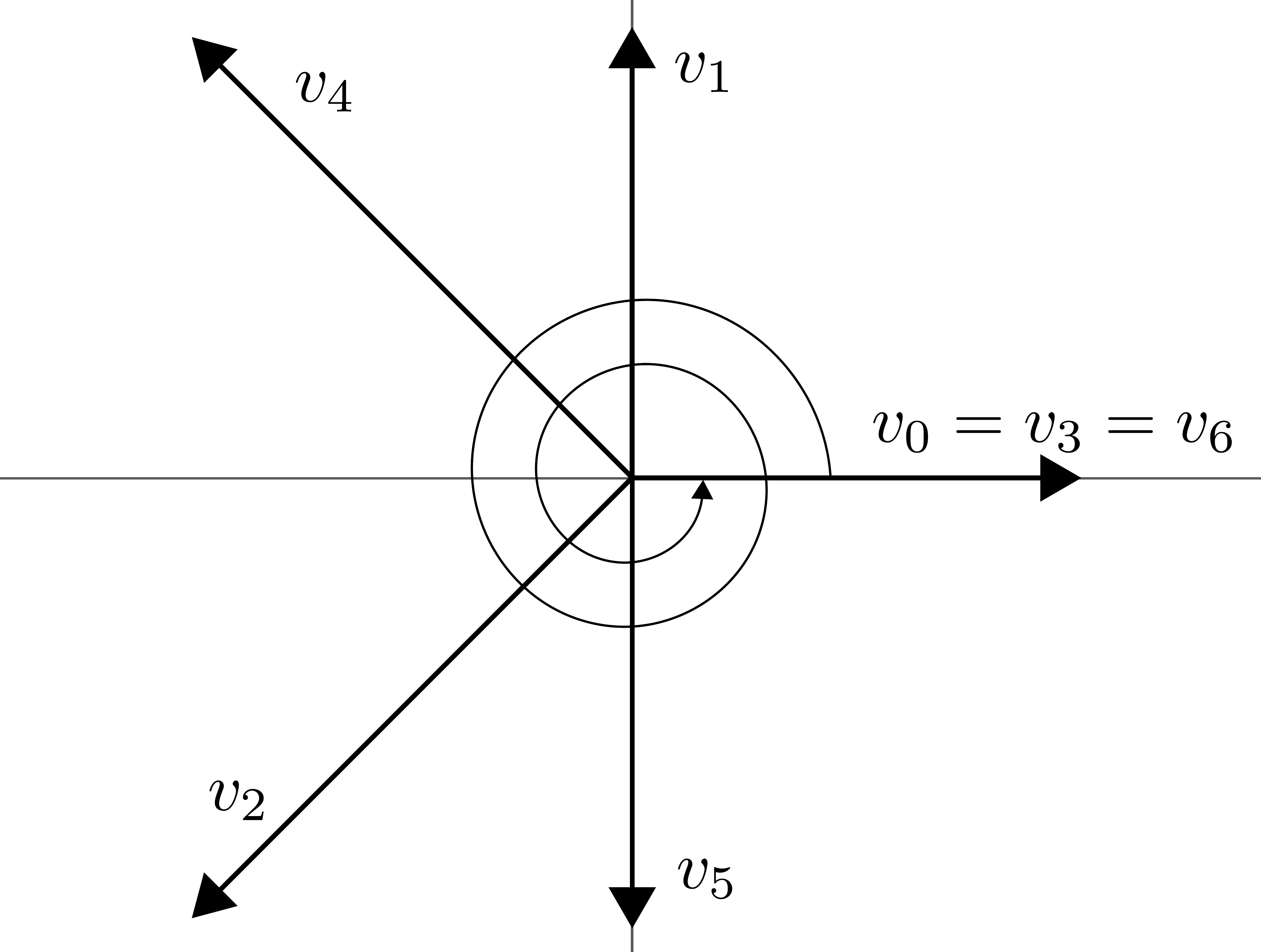}
 \caption{These vectors do not form a fan because they are not labeled in counter-clockwise order.} \label{fig_windtwice}
 \end{figure}

\begin{Example}

 Consider the sequence of integers $a_0=-1$, $a_1=-1$, $a_2=-2$, $a_3=-1$, $a_4=-1$, $a_5=0$ and notice that
 \begin{gather*}
 ST^{-1}ST^{-1}ST^{-2}ST^{-1}ST^{-1}ST^{0}=_{\sltz}I,
 \end{gather*}
 so equation~\eqref{eqn_st} is satisf\/ied but these integers do not correspond to any toric fan. This can be seen by considering the vectors:
 \begin{gather*}
 v_0 = \vect{1}{0}, \qquad\!\! v_1 = \vect{0}{1},\qquad\!\! v_2 = \vect{-1}{-1},\qquad\!\!
 v_3 = \vect{1}{0},\qquad\!\! v_4 = \vect{-1}{1}, \qquad\!\! v_5 = \vect{0}{-1}
 \end{gather*}
 which travel twice around the origin (a formal def\/inition of this is given in Def\/inition \ref{def_circleorigin}), see Fig.~\ref{fig_windtwice}, and form one possible sequence of vectors which can be obtained from the given integers $a_0, \ldots, a_5$. Since the action of $\sltz$ preserves the property of being labeled in counter-clockwise order this means that none of the sequences of vectors which can be obtained from these integers is a toric fan.
\end{Example}

So we need extra information that is not captured by viewing this word in $S$ and $T$ as an element of $\sltz$. For a more obvious example notice that even though they are equal in $\sltz$ we can see that $S^4$ corresponds to a toric fan while $S^8$ does not. From~\cite[p.~44]{Fultontoric} we know that integers $(a_0, \ldots, a_d) \in \Z^d$ which satisfy equation~\eqref{eqn_intromatr} correspond to a toric fan if and only if
\begin{gather*} a_0 + \cdots + a_{d-1} = 3d - 12,\end{gather*} so we would like to prove that
\begin{gather*}
 \frac{3d - \sum\limits_{i=0}^{d-1}{a_i}}{12}
\end{gather*}
is the number of times that the vectors corresponding to $(a_0, \ldots, a_{d-1})$ circle the origin. In order to prove this we will need some more terminology, and in order to keep track of the extra information about circling the origin we will need to lift to a group which projects onto $\sltz$.

Consider instead the free group with generators $S$ and $T$. This group is isomorphic to the free product of $\Z$ with itself so we will denote
it by $\Z*\Z$. We know $\sltz$ is a quotient of $\Z * \Z$ by Lemma~\ref{lem_sltzpresn} so there exists a natural projection map $\pi_1\colon \Z * \Z \to \sltz$. Also def\/ine a map $w\colon \Z *\Z\to\Z$ to be the unique homomorphism such that
\begin{gather*}
 w(S_{\Z*\Z}) = 3\qquad \text{and} \qquad w(T_{\Z*\Z}) = -1.
\end{gather*}
Given any word in the kernel of the projection from $\Z*\Z$ to $\sltz$ we will see that $w$ evaluates to 12 times the image of the map to $\pi_1(\sltr)$ that is obtained by lifting the word in the natural way to the universal cover of $\sltr$, and $w$ is the unique homomorphism on all of $\Z*\Z$ with this property. Given a toric fan with associated integers $(a_0, \ldots, a_{d-1})\in\Z^d$ we will show that
\begin{gather*}w\big(\big(ST^{a_0}\cdots ST^{a_{d-1}}\big)_{\Z*\Z}\big)=12.\end{gather*}
Both $\pi_1$ and $w$ factor over the same group $G$ which is the f\/iber product of~$\sltz$ and $\Z$ over~$\Z/(12)$. Now we can see that we wanted the particular presentation of~$\sltz$ from Lemma~\ref{lem_sltzpresn} so that the relationship between~$G$ and~$\sltz$ would be clear. This discussion is made precise in the following proposition.

\begin{Proposition}\label{prop_diagram}
The following diagram commutes:
\begin{gather*}
\xymatrix{
 &\Z*\Z \ar[dl]_{\pi_1} \ar[d]_{\pi_2} \ar@/^1pc/[dd]^{w} \\
\sltz \ar[d]_{w_{\sltz}} & G \ar[l]_-{\pi_3} \ar[d]_{w_G} \\
\Z/(12) &\Z. \ar[l]^-{\pi_4} }
\end{gather*}
The group $G$ is the fiber product of $\sltz$ and $\Z$ over $\Z/(12)$ and satisfies
\begin{gather*}
 G \cong \big\langle S,T \,|\, STS = T^{-1}ST^{-1} \big\rangle,
\end{gather*}
each of $\pi_1$, $\pi_2$, $\pi_3$, and $\pi_4$ is a projection, and $w\colon \Z*\Z\to\Z$, $w_G\colon G\to\Z$, $w_{\sltz}\colon \sltz\to\Z/(12)$ are given by the same formal expression
\begin{gather}\label{eqn_formalexpr}
 S^{b_0}T^{a_0}\cdots S^{b_\ell}T^{a_\ell} \mapsto 3\sum_{i=0}^\ell b_i - \sum_{i=0}^\ell a_i.
\end{gather}
\end{Proposition}

\begin{proof} It can be seen that the map $w_{\sltz}\colon \sltz \to \Z/(12)$ is well-def\/ined by noting that both relations in $\sltz$ as presented in Lemma~\ref{lem_sltzpresn} preserve the value of the formula (equation~\eqref{eqn_formalexpr}) up to a multiple of 12. Similarly, since the relation $STS=_{\sltz}T^{-1}ST^{-1}$ preserves the value of the equation we know that $w_G$ is well-def\/ined. Since each of these functions to~$\Z$ or~$\Z/(12)$ is given by the same formal expression and since each $\pi$ is a quotient map, the diagram commutes.

To show that $G$ with the associated maps is isomorphic to the f\/iber product of~$\sltz$ and~$\Z$ over $\Z/(12)$ we must only show that $w_G$ restricted to the f\/ibers is bijective. That is, we must show that
\begin{gather*}
 w_G\restriction_{\pi_3^{-1}(A)}\colon \ \pi_3^{-1}(A)\to \pi_4^{-1}(w_{\sltz}(A))
\end{gather*}
is a bijection for each $A\in\sltz$. To show it is surjective, f\/ix $A\in\sltz$ and let $g\in\pi_3^{-1}(A)$. For any $k\in\Z$ notice $\pi_3\big(S^{4k}g\big) = \pi_3(g)$ so $S^{4k}g\in\pi_3^{-1}(A)$. Since $w_G\big(S^{4k}g\big) = w_G(g) + 12k$ and $\pi_4^{-1}(w_{\sltz}(A))=w_G(g)+12\Z$ the map is surjective. To show it is injective it is suf\/f\/icient to consider only $A = I$. Since $\big(S^4\big)_G$ is in the center of $G$ we know that $\sltz \cong G/\big(S^4\big)$ so $\pi_3^{-1}(I)=\big\{ S^{4k}_G \,|\, k\in\Z \big\}$. Since $w_G\big(S^{4k}_G\big) = 12k$ we know for each choice of $k$ this maps to a~distinct element of~$\Z$.
\end{proof}

There is another useful sense in which $G$ is an unwinding of $\sltz$. While $\sltz$ is discrete, and thus does not have a natural cover, it sits inside the group $\mathrm{SL}_2(\R)$, which has a universal cover. We claim that~$G$ is the preimage of $\sltz$ inside of the universal cover of~$\mathrm{SL}_2(\R)$.
\begin{Proposition}\label{prop_universalcover}
The group $G$ is isomorphic to the preimage of $\sltz$ within the universal cover of $\mathrm{SL}_2(\R)$.
\end{Proposition}
\begin{proof} Let $G'$ be the preimage of $\sltz$ in the universal cover of~$\mathrm{SL}_2(\R)$. We note that there exists a homomorphism,
 $\phi$ from $G$ to $G'$ def\/ined by
\begin{gather*}
\phi(S_G) = \left(\begin{matrix} \cos\left(\frac{\pi t}{2}\right) & -\sin\left(\frac{\pi t}{2}\right)\vspace{1mm}\\
\sin\left(\frac{\pi t}{2}\right) & \cos\left(\frac{\pi t}{2}\right)\end{matrix}\right)_{0\leq t \leq 1}, \qquad \phi(T_G) = \matr{1}{t}{0}{1}_{0\leq t \leq 1},
\end{gather*}
where the paths given are to represent elements of the universal cover of $\mathrm{SL}_2(\R)$. Consider how $\phi(S_G)\phi(T_G)\phi(S_G)$ and $\phi(T_G)^{-1}\phi(S_G)\phi(T_G)^{-1}$ act on any basis of~$\Z^2$. Both paths will take the basis to the same point without traversing completely around the origin (i.e., there is a ray starting at the origin in~$\R^2$ that the basis does not cross as either path acts on it). From this fact it is easy to see that they are homotopic relative to endpoints, and thus~$\phi$ does actually def\/ine a homomorphism. We have left to show that $\phi$ def\/ines an isomorphism. To show this, we note that each of $G$ and $G'$ have obvious surjections~$\pi$ and~$\pi'$ to~$\sltz$. Furthermore, $\pi = \pi'\circ \phi$. Thus, to show that $\phi$ is an isomorphism, it suf\/f\/ices to show that $\phi\colon \ker(\pi)\rightarrow\ker(\pi')$ is an isomorphism.

Notice that $\ker(\pi)$ is $\langle S_G^4 \rangle$. On the other hand, $\ker(\pi')\cong \pi_1(\mathrm{SL}_2(\R))=\Z$, and is generated by $\phi(S_G)^4.$ This completes the proof.
\end{proof}

We will see in Corollary~\ref{cor_toricfaneqn} that (as shown in \cite[equation~(3)]{PR2000}) there is a one to one correspondence between toric fans up to the action of $\sltz$ and lists of integers $a_0, \ldots, a_{d-1}\in\Z$ satisfying
\begin{gather}\label{eqn_stfan}
 ST^{a_0}\cdots ST^{a_d-1} =_G S^4.
\end{gather}
Equation \eqref{eqn_stfan} is a ref\/inement of equation~\eqref{eqn_st} which implies both that the successive pairs of vectors form a basis of~$\Z^2$ and that the vectors are labeled in counter-clockwise order. In Proposition~\ref{prop_standardformstfan} we produce an analogous equation for semitoric fans.

Now we would like to simplify these toric fans. We will understand which integers $a_0, \ldots, a_{d-1}$ are possible in an element $(ST^{a_0}\cdots ST^{a_{d-1}})_G\in G$ corresponding to a toric fan by stu\-dying $\psl\cong\sltz/(-I)$. We again use Notation~\ref{notationST} and write
$S_{\psl}$ and $T_{\psl}$ as the elements in $\psl$ corresponding to $S_\sltz$ and $T_\sltz$ in $\sltz$. The following lemma is important for this and will also be useful later on when classifying semitoric fans. If $(ST^{a_0}\cdots ST^{a_{d-1}})_G$ $\in G$ projects to the identity in $\sltz$ then by Lemma \ref{lem_psltz} we can see when one of the exponents must be in the set $\{-1, 0, 1\}$. In any of these cases, we will be able to use relations in $G$ to help simplify the expression.

\begin{Lemma}\label{lem_psltz}
 Suppose that \begin{gather}\label{eqn_lempsltz} ST^{a_0}\cdots ST^{a_{d-1}}=_{\psl} I\end{gather} for some $d\in\Z$, $d>0$. Then if $d\geq3$ there exist $i,j\in\Z$ satisfying $0\leq i < j \leq d-1$ such that $a_i, a_j \in \{-1, 0, 1\}$. Furthermore:
 \begin{enumerate}\itemsep=0pt
 \item[$1)$]\label{part_psltzone}
 if $d>3$ then $i,j$ can be chosen such that $i \neq j-1$ and $(i,j)\neq(0,d-1)$;
 \item[$2)$] if $d=3$ then $a_0=a_1=a_2=1$ or $a_0=a_1=a_2=-1$;
 \item[$3)$] if $d=2$ then $a_0=a_1=0$;
 \item[$4)$] if $d=1$ then equation \eqref{eqn_lempsltz} cannot hold.
 \end{enumerate}
\end{Lemma}

Note that part (1) 
is the statement that $i$ and $j$ are not consecutive in the cyclic group $\Z/(d)$. Of course, such $i$ and $j$ may be chosen if three or more elements of the list $a_0, \ldots, a_{d-1}$ are in the set $\{-1,0,1\}$.

\begin{proof} It is well known that $\psl$ acts faithfully on the real projective line $\R\mathrm{P}^1 = \R\cup\{\infty\}$ by linear fractional transformations: \begin{gather*}\matr{a}{b}{c}{d} (x) = \frac{ax+b}{cx+d} \quad \text{for $x\in\R$} \qquad \text{and} \qquad \matr{a}{b}{c}{d} (\infty) = \frac{a}{c}.\end{gather*}

 Let $d>3$. To f\/ind a contradiction suppose that at most two of $a_0, \ldots, a_{d-1}$ are in $\{-1, 0, 1\}$, and also suppose that if
 there are two in $\{-1, 0, 1\}$ that they are consecutive or indexed by $0$ and $d-1$. Notice that \begin{gather*}ST^{a_0}\cdots ST^{a_{d-1}}=_{\psl}I \qquad \text{implies that} \qquad ST^{a_1}\cdots ST^{a_{d-1}}ST^{a_0} =_{\psl} I\end{gather*} by conjugating each side with $(ST^{a_0})_{\psl}$. This conjugation method and renumbering the integers can be used to assure that $a_i \notin \{-1, 0, 1\}$ for $i=1, \ldots, d-2$. Since this expression is equal to the identity in $\psl$ it acts trivially on $\R\cup \{\infty\}$. In particular, we have
 \begin{gather*}\big(ST^{a_0}\cdots ST^{a_{d-1}}\big)_{\psl} (\infty) = \infty.\end{gather*}
Notice that $S_{\psl}(x) = \nicefrac{-1}{x}$ and $T_{\psl}^a(x) = x+a$ for $a\in\Z$. Further notice that for any $a\in\Z\setminus\{-1,0,1\}$ and $x\in(-1,1)\setminus\{0\}$ we have $(ST^a)_{\psl} (x) \in(-1,1)\setminus\{0\}$. We see that $(ST^{a_{d-1}})_{\psl}(\infty) = 0$ and since $a_{d-2}\notin \{-1, 0, 1\}$ we know $(ST^{a_{d-2}})_{\psl}(0) \in (-1,1)\setminus\{0\}$. Putting these facts together we have
 \begin{align*}
 \big(ST^{a_0}\cdots ST^{a_{d-1}}\big)_{\psl} (\infty) &= \big(ST^{a_0}\cdots ST^{a_{d-2}}\big)_{\psl} (0) &\\
 &=\big(ST^{a_0}\cdots ST^{a_{d-3}}\big)_{\psl}(x) &\text{for some }x\in (-1,1)\setminus \{0\}\\
 &=\big(ST^{a_0}\big)_{\psl} (y) &\text{for some } y\in (-1,1)\setminus\{0\}\\
 &=\frac{-1}{y+a_0}\neq \infty.
 \end{align*}
 This contradiction f\/inishes the $d > 3$ case.

 If $d=3$ essentially the same result holds except that it is not possible to choose two elements that are non-consecutive. Notice
 \begin{gather*}\big(ST^{a_0}ST^{a_1}ST^{a_2}\big)_{\psl}(\infty) = \big(ST^{a_0}ST^{a_1}\big)_{\psl} (0) = -\left(\frac{-1}{a_1} + a_0\right)^{-1}.\end{gather*}
 For this function to be the identity we would need
 \begin{gather*}\frac{-1}{a_1} +a_0=0.\end{gather*}
 This implies that $a_0 a_1 = 1$ so since they are both integers we have $a_0 = a_1 = \ep$ where $\ep\in\{-1,1\}$. Conjugating by $(ST^{a_0})_{\psl}$, we f\/ind symmetrically, that $a_1a_2=1$, and thus that $a_0=a_1=a_2=\pm 1$.

 If $d=2$ then we have
 \begin{gather*}\big(ST^{a_0}ST^{a_1}\big)_{\psl}(\infty) = \big(ST^{a_0}\big)_{\psl}(0)=\frac{-1}{a_0},\end{gather*}
 so we must have $a_0 = 0$ and then
 \begin{gather*} \big(ST^0ST^{a_1}\big)_{\psl}(x) =\big(S^2 T^{a_1}\big)_{\psl}(x)=x+a_1,\end{gather*}
 so we are also forced to have that $a_1=0$, as stated in the lemma. If $d=1$ then $(ST^{a_0})_{\psl}(\infty)$ $ = 0$ for any choice of $a_0 \in \Z$, so there are no solutions.
\end{proof}

Note that $ST^0ST^0\neq_{\sltz}I$ and $ST^1ST^1ST^1\neq_{\sltz}I$, so we have the following corollary.

\begin{Corollary}\label{cor_dbiggerthree}
 If $d\in\Z$, $d>0$, and $a_0, \ldots, a_{d-1}\in\Z$ are such that
 \begin{gather*}
 ST^{a_0}\cdots ST^{a_{d-1}}=_{\sltz}I,
 \end{gather*}
 then $d\geq 3$.
 Furthermore, if $d=3$ then $a_0=a_1=a_2=-1$.
\end{Corollary}

\begin{Remark}
 In terms of toric fans, Corollary~\ref{cor_dbiggerthree}
 is saying that every toric fan has at least three vectors and,
 up to the action of $\sltz$, there is only one toric fan
 with exactly three vectors (the fan of $\mathbb{CP}^2$).
\end{Remark}

\section{Toric fans}\label{sec_toric}

Let $a_0, a_1, \ldots, a_{d-1}\in\Z$ be a collection of integers such that
\begin{gather*}
 ST^{a_0}\cdots ST^{a_{d-1}} =_{\sltz} I.
\end{gather*}
This means that
\begin{gather*}
 ST^{a_0}\cdots ST^{a_{d-1}} =_G S^{4k}\qquad \text{for some} \ k\in\Z
\end{gather*}
by Proposition \ref{prop_diagram}. We claim that these integers correspond to a toric fan if and only if $k=1$.

The idea is that from $a_0, a_1, \ldots, a_{d-1}$ we can always construct a sequence of vectors for which every pair of consecutive vectors is a basis of $\Z^2$, but $k=1$ precisely when these vectors are labeled in counter-clockwise order. The only relation in the group~$G$, which is $STS =_G T^{-1}ST^{-1}$, preserves the number of times the vectors circle the origin, so it is natural to study winding in $G$. Now we make this idea precise.

\begin{Lemma}
 Let $g\in\ker (\pi_3)$. Then $\frac{w_G (g)}{12} \in \Z$.
\end{Lemma}

\begin{proof}
 Since $\pi_3 (g) = I $ we know $w_{\sltz} \circ \pi_3 (g) = 0$, so
 by Proposition \ref{prop_diagram} $\pi_4 \circ w_G (g) = 0$. Thus
 $w_G (g)\in \ker (\pi_4) = \{ 12k \,|\, k\in\Z \}.$
\end{proof}

Recall that $G$ is isomorphic to the preimage of $\sltz$ in the universal cover of $\mathrm{SL}_2 (\R)$ by Proposition~\ref{prop_universalcover}. Let $\phi$ be the isomorphism from $G$ to its image in the universal cover of $\mathrm{SL}_2(\R)$ with
\begin{gather*}
 \phi(S_G) = \left(\begin{matrix} \cos\left(\frac{\pi t}{2}\right) & -\sin\left(\frac{\pi t}{2}\right)\vspace{1mm}\\
 \sin\left(\frac{\pi t}{2}\right)& \cos\left(\frac{\pi t}{2}\right)\end{matrix}\right)_{0\leq t \leq 1}, \qquad \phi(T_G) = \matr{1}{t}{0}{1}_{0\leq t \leq 1}.
\end{gather*}
This means to each element of the kernel of $\pi_3$ we can associate a closed loop based at the identity in $\mathrm{SL}_2 (\R)$ denoted~$\phi(g)$. The fundamental group $\pi_1(\sltr)$ is isomorphic to~$\Z$ and is generated as $\big\langle \phi\big(S^4_G\big) \big\rangle$, so let $\psi\colon \pi_1 (\sltr) \to \Z$ be the isomorphism with $\psi (\phi (S^4_G)) = 1$.

\begin{Lemma}\label{lem_windingnumber}
 Let $g\in\ker (\pi_3\colon G \to \sltz)$. Then \begin{gather*}\psi\circ\phi (g)= \frac{w_G(g)}{12}.\end{gather*}
\end{Lemma}

\begin{proof} Since $\ker(\pi_3)$ is generated by $S^4_G$, it suf\/f\/ices to check that $\psi\big(\phi\big(S^4_G\big)\big)=\frac{w_G(S^4_G)}{12}=1$, but this holds by def\/inition.
\end{proof}

\begin{Definition}\label{def_windingnumber} Def\/ine $W\colon \ker (\pi_3) \to \Z$ by
 \begin{gather*}
 W(g) = \frac{w_G (g)}{12}.
 \end{gather*} We call $W(g)$ the \emph{winding number} of $g\in\ker (\pi_3)$.
\end{Definition}

\begin{Definition}\label{def_circleorigin} Let \begin{gather*}(v_0=v_d, v_1, \ldots, v_{d-1})\in\big(\Z^2\big)^d\end{gather*} with $\det(v_i, v_{i+1})>0$
 for $i=0\ldots, d-1$. We def\/ine \emph{the number of times $(v_0, \ldots, v_{d-1})$ circles the origin} to be the winding number of the piecewise linear path in $\big(\R^2\big)^* = \R^2\setminus\{(0,0)\}$ produced by concatenating the linear paths between $v_i$ and $v_{i+1}$ for $i=0,\ldots, d-1$.
\end{Definition}

\begin{Lemma}\label{lem_toricwinding} Let $a_0, \ldots, a_{d-1}\in\Z$ such that $ST^{a_0}\cdots ST^{a_{d-1}} =_{\sltz} I$ and let $v_0, v_1\in\Z^2$
 such that $\det(v_0, v_1) = 1$. Then $d\geq 3$ so we may define $v_2, \ldots v_{d-1}$ by
 \begin{gather*}
 v_{i+2} = -v_i + a_i v_{i+1},
 \end{gather*}
 where $v_d = v_0$ and $v_{d+1} = v_1$. The winding number $W((ST^{a_0}\cdots ST^{a_{d-1}})_G)\in\Z$ is the number of times that $(v_0, \ldots, v_{d-1})\in\big(\Z^2\big)^d$ circles the origin.
\end{Lemma}

\begin{proof} Corollary~\ref{cor_dbiggerthree} states that $d\geq 3$ in this situation. For $i=1, \ldots, d-2$ let $A_i$ be the matrix $[v_i,v_{i+1}]$, let $A_{d-1} = [v_{d-1},v_0]$, and recall that $A_{i+1}=A_i(ST^{a_i})_\sltz$. Identify~$G$ as a~subgroup of the universal cover of~$\mathrm{SL}_2(\R)$ by choosing the identity as a basepoint. Thus, for each $i$ there is a natural path between $I$ and $(ST^{a_i})_\sltz$, and multiplying this path on the left by $A_i$ we obtain a path from $A_i$ to $A_i (ST^{a_i})_\sltz = A_{i+1}$. These paths for $i=0, \ldots, d-1$ can be concatenated to form a path from $A_0$ to $A_0$. Projecting this path into the f\/irst column vector of the appropriate matrix gives a path in~$\big(\R^2\big)^*$. We claim that this path is homotopic to the path formed by taking line segments between $v_i$ and $v_{i+1}$. This is easily verif\/ied because both paths between $v_i$ and $v_{i+1}$ travel counterclockwise less than a full rotation.

We know that $W((ST^{a_0}\cdots ST^{a_{d-1}})_G)$ equals $\psi$ of the path in $\mathrm{SL}_2(\R)$ by Lemma~\ref{lem_windingnumber}, and now we need to show that this equals the winding number in $\big(\R^2\big)^*$ of the f\/irst column vectors. This holds because projection onto the f\/irst column vector of a matrix induces an isomorphism of fundamental groups between $\mathrm{SL}_2(\R)$ and $\big(\R^2\big)^*$ which sends the generator $[\phi(S_G)^4]$ to the homotopy class of a loop which circles counter-clockwise once around the origin, a path with winding number~1.
 \end{proof}

\begin{Corollary}\label{cor_toricfaneqn}
 There exists a bijection from the set of all sequences $(a_0, \ldots, a_{d-1})\in\Z^d$, $d>0$, which satisfy \begin{gather*} ST^{a_0}\cdots ST^{a_{d-1}} =_G S^4\end{gather*} to the collection of all toric fans modulo the action of $\sltz$. This bijection sends $(a_0, \ldots, a_d)\in\Z^d$ to the equivalence class of fans \begin{gather*}\big\{(v_0=v_d, v_1 = v_{d+1}, v_2, \ldots, v_{d-1})\in \big(\Z^2\big)^d\,|\, v_0, v_1\in\Z^2, \, \det(v_0, v_1)=1\big\}\end{gather*} in which \begin{gather*}v_{i+2} = -v_i + a_i v_{i+1}\end{gather*} for $i=0, \ldots, d-1$.
\end{Corollary}

\begin{proof} Let $(v_0 = v_d, v_1, \ldots, v_{d-1}) \in\big(\Z^2\big)^d$ be a toric fan. That is, $\det(v_i, v_{i+1})=1$ for each $i=0, \ldots, d-1$ and the vectors are labeled in counter-clockwise order. It is shown in Section~\ref{sec_algsetup}, equation~\eqref{eqn_st} that associated integers $(a_0, \ldots, a_{d-1})\in(\Z)^d$ exist such that
 \begin{gather*}
 ST^{a_0}\cdots ST^{a_{d-1}} =_{\sltz} I,
 \end{gather*}
 which means
 \begin{gather*}
 ST^{a_0}\cdots ST^{a_{d-1}} =_G S^{4k}
 \end{gather*}
 for some $k\in\Z$ with $k\geq0$. By Lemma \ref{lem_toricwinding} we know
 \begin{gather*}
 W\big(\big(ST^{a_0}\cdots ST^{a_{d-1}}\big)_G\big)=1,
 \end{gather*}
so that the vectors will be labeled in the correct order for it to be a~fan. Thus, $W\big(S^{4k}_G\big) = 1$ but $W\big(S^{4k}_G\big) = k$ so $k=1$. Notice such a construction is well-def\/ined on equivalence classes of toric fans because the integers are prescribed via linear equations and fans in a~common equivalence class are related by a linear map.

 Now suppose that $(a_0, \ldots, a_{d-1}) \in \Z^d$ satisfy $ST^{a_0}\cdots ST^{a_{d-1}} =_G S^4$ and def\/ine $(v_0, \ldots, v_{d-1})$ $\in\big(\Z^2\big)^d$ by
 \begin{gather*}
 v_{i+2} = -v_i + a_i v_{i+1},
 \end{gather*}
 where $v_0, v_1\in\Z^2$ are any two vectors for which $\det(v_0, v_1)=1$. Then for each $i=0, \ldots, d-1$ we have
 \begin{gather*}
 \det(v_{i+1}, v_{i+2}) = \det\matr{0}{-1}{1}{a_i}\det(v_i,v_{i+1}) = \det(v_i,v_{i+1}),
 \end{gather*}
 so by induction all of these determinants are 1. By Lemma \ref{lem_toricwinding}, the path connecting adjacent vectors wraps around the origin only once, and since each $v_{i+1}$ is located counterclockwise of $v_i$, we have that the $v_i$'s must be sorted in counterclockwise order.

 It is straightforward to see that these constructions are inverses of one another.
\end{proof}

Now that we have set up the algebraic framework the following results are straightforward to prove. First we prove that any fan with more than four vectors can be reduced to a fan with fewer vectors. This result is well-known but the following proof is new.

\begin{Lemma}[{\cite[Claim, p.~43]{Fultontoric}}]\label{lem_fulton}
 If $(v_0=v_d, \ldots, v_{d-1})\in\big(\Z^2\big)^d$ is a toric fan with $d>4$ then there exists some $i\in \{0, \ldots, d-1\}$ such that $v_i = v_{i-1} + v_{i+1}$.
\end{Lemma}

\begin{proof} By Corollary \ref{cor_toricfaneqn} we know that to the fan $(v_0=v_d, \ldots, v_{d-1})$ there is an associated list of integers ${a_0, \ldots, a_{d-1}\in\Z}$ such that $v_{i+2} = -v_i + a_i v_{i+1}$ and
 \begin{gather}\label{eqn_lemfulton}
 ST^{a_0}\cdots ST^{a_{d-1}} =_G S^4.
 \end{gather}
 We must only show that for some $i\in\Z$ we have $a_i = 1$. Since $S^4 =_{\psl} I$ we can use Lemma~\ref{lem_psltz} to conclude that there exist $i,j\in\Z$ satisfying $0\leq i<j-1\leq d-2$ such that $a_i, a_j\in\{-1, 0, 1\}$ and $(i,j)\neq (0, d-1)$. By way of contradiction assume that $a_i, a_j \in\{-1, 0\}$. Conjugate equation \eqref{eqn_lemfulton} by $ST^{a_n}$ for varying $n\in\Z$ to assure that $i\neq 0$ and $j\neq d-1$. Then at each of these values we may use either $ST^0S =_G S^2$ or $ST^{-1}S =_G S^2 TST$ to reduce the number of $ST$-pairs by one or two and produce a factor of $S^2_G$, which we can move to the front of the word because $S^2_G$ is in the center of $G$. These reductions do not interfere with one another because the values in question are not adjacent. So we end up with
 \begin{gather*}
 S^4 ST^{b_0}\cdots ST^{b_{\ell-1}} =_G S^4
 \end{gather*}
 and thus
 \begin{gather*}
 ST^{b_0}\cdots ST^{b_{\ell-1}} =_G I,
 \end{gather*}
where $\ell \geq 3$ by Corollary~\ref{cor_dbiggerthree}. This implies that $W((ST^{b_0}\cdots ST^{b_{\ell-1}})_G)=0$ and thus, by Lemma~\ref{lem_toricwinding}, that any corresponding collection of vectors winds no times about the origin. However, this is impossible since for such a sequence of vectors each vector is always counterclockwise from the previous vector and $\ell>1$.
\end{proof}

The case in which a vector in the fan is the sum of the adjacent vectors is important because this means the fan is the result of corner chopping a fan with fewer vectors in it. Now that we have the proper algebraic tools, we will be clear about the specif\/ics of the corner chopping and reverse corner chopping operations.

Suppose $(v_0 = v_d, v_2, \ldots, v_{d-1}) \in \big(\Z^2\big)^d$ is a toric fan with associated integers $(a_0, \ldots, a_{d-1}) \in\Z^d$. Then \begin{gather*} v_{i+2} = -v_i + a_{i}v_{i+1}\end{gather*} so if $a_i = 1$ then we have that $v_{i+1} = v_i + v_{i+2}$. Now we see that in this case \begin{gather*} \det(v_i, v_{i+2}) = \det(v_i,-v_i )+\det(v_i,v_{i+1}) = 1,\end{gather*} so \begin{gather*} (w_0 = v_0, \ldots, w_i = v_i, w_{i+1} = v_{i+2}, \ldots, w_{i-2}=v_{d-1})\in\big(\Z^2\big)^{d-1}\end{gather*} is also a fan. Next notice \begin{gather*} -w_{i}+(a_{i+1}-1)w_{i+1} = -(v_i+v_{i+2}) + a_{i+1} v_{i+2} = -v_{i+1} + a_{i+1} v_{i+2} = w_{i+2}\end{gather*} and
\begin{gather*} -w_{i-1}+(a_{i-1}-1)w_{i} = (-v_{i-1}+a_{i-1} v_{i}) - v_i = v_{i+1} - v_i = w_{i+1},\end{gather*}
so this new fan has associated to it the tuple of integers \begin{gather*}(a_0, \ldots, a_{i-1}-1, a_{i+1}-1, \ldots, a_{d-1})\in \Z^{d-1}.\end{gather*} An occurrence of $1$ from the original tuple of integers has been removed and the adjacent integers have been reduced by $1$. Algebraically, this move corresponds to the relation $STS=_G T^{-1}ST^{-1}$. Geometrically this move corresponds to the inverse of chopping a corner from the associated polygon (as is shown in Fig.~\ref{fig_minimalmodels}) and from the point of view of symplectic geometry this corresponds to an equivariant symplectic blowup. The corner chopping of a toric polygon is done such that the new face of the polygon produced has inwards pointing normal vector given by the sum of the adjacent inwards pointing primitive integer normal vectors.

Now we can see that Lemma \ref{lem_fulton} tells us that fans with f\/ive or more vectors are the result of corner chopping a fan with fewer vectors. We will next classify all possible fans with fewer than f\/ive vectors.

\begin{Lemma}\label{lem_minmodels}
 Suppose that integers $a_0, \ldots, a_{d-1}\in\Z$ satisfy \begin{gather}\label{eqn_lemminmodels}ST^{a_0}\cdots ST^{a_{d-1}} =_G S^4\end{gather} for some $d\in\Z$, $d\geq0$.
 \begin{enumerate}\itemsep=0pt
 \item[$1.$] If $d=4$ then up to a cyclic reordering the set of integer quadruples which satisfy this equation is exactly $a_0 = 0$, $a_1= k$, $a_2 = 0, a_3 = -k$ for each $k\in\Z$.
 \item[$2.$] 
 If $d=3$ then $a_0 = a_1 = a_2 = -1$.
 \item[$3.$] 
 If $d<3$ then there do not exist integers satisfying equation \eqref{eqn_lemminmodels}.
 \end{enumerate}
\end{Lemma}
\begin{proof} Notice $ST^{a_0}\cdots ST^{a_{d-1}} =_G S^4$ implies that $ ST^{a_0}\cdots ST^{a_{d-1}} =_{\sltz} I$ so Corollary~\ref{cor_dbiggerthree} implies items~(2) and~(3). Moreover, we also have that
\begin{gather*}
 ST^{a_0}\cdots ST^{a_{d-1}} =_{\psl} I,
\end{gather*}
so we can apply Lemma~\ref{lem_psltz}. Suppose that $d=4$. Lemma \ref{lem_psltz} tells us that at least one of the $a_i$ is in the set $\{ -1, 0, 1\}$. By conjugation (which cyclically permutes the order of the integers) we may assume that $a_0 \in \{ -1, 0, 1 \}$. If $a_0 = 1$ then \begin{gather*} STST^{a_1}ST^{a_2}ST^{a_3} =_G ST^{a_1-1}ST^{a_2}ST^{a_3-1},\end{gather*} so for this to equal $S^4_G$ in $G$ we must have $a_1-1 = a_2 = a_3-1 = -1$ by the case of $d=3$. It is straightforward to check that $STSST^{-1}S =_G S^4$ so we have found the required solution.

 If $a_0 = -1$ then notice
\begin{gather*}
 ST^{-1}ST^{a_1}ST^{a_2}ST^{a_3}=_G S^4 \qquad \text{implies} \qquad ST^{a_1 + 1}ST^{a_2}ST^{a_3+1}=_G S^2 =_{\psl} I,
\end{gather*} so by Lemma~\ref{lem_psltz} we must have $a_1+1=a_2=a_3+1=\pm 1.$ This time, if $a_1+1=a_2=a_3+1= -1$ then equation~\eqref{eqn_lemminmodels} does not hold, since the left side will equal $S^6$, but if $a_1+1=a_2=a_3+1=1$ then the equation holds. So we have found another solution, $ST^{-1}SSTS=_G S^4$, which has the form described in the statement of the lemma.

 Finally, suppose that $a_0 = 0$. Notice
 \begin{gather*}
 ST^0ST^{a_1}ST^{a_2}ST^{a_3}=_G S^4 \qquad \text{implies} \qquad ST^{a_2}ST^{a_1+a_3} =_G S^2 =_{\psl} I,
 \end{gather*}
 so we can use Lemma \ref{lem_psltz} to conclude that we need $a_2 = a_1 + a_3 = 0$. Let $a_1 = k\in\Z$. Now we have that \begin{gather*} ST^0ST^kST^0ST^{-k}=_G S^4\end{gather*} for any $k\in\Z$. Finally, observe that the other two possibilities we derived in the $d=4$ case are just reorderings of this one with $k=1$.
\end{proof}

\begin{Definition}\label{def_trisquaretrap} A \emph{Delzant triangle} is the convex hull of the points $(0,0)$, $(0, \la)$, $(\la, 0)$ in $\R^2$ for any $\la>0$. A \emph{Hirzebruch trapezoid} with parameter $k\in\Z_{\geq0}$ is the convex hull of $(0,0)$, $(0,a)$, $(b, a)$, and $(b+ak,0)$ in~$\R^2$ where $a,b>0$. A Hirzebruch trapezoid with parameter zero is a~\emph{rectangle}.
\end{Definition}

These are shown in Fig.~\ref{fig_minimalmodels}. So we see that the fan corresponding to any Delzant triangle is \begin{gather*} \left(\vect{1}{0}, \vect{0}{1}, \vect{-1}{-1}\right)\end{gather*} with associated integers $(-1,-1,-1)$ and the fan corresponding to a Hirzebruch trapezoid with parameter $k$ is \begin{gather*}\left( \vect{0}{1}, \vect{-1}{-k}, \vect{0}{-1}, \vect{1}{0}\right)\end{gather*} with associated integers $(0, k, 0, -k)$. The following Theorem is immediate from Lemmas~\ref{lem_fulton} and~\ref{lem_minmodels}.

\begin{Theorem}[\cite{Fultontoric}]\label{thm_toricpoly}
 Every Delzant polygon can be obtained from a polygon $\sltz$-equivalent to a Delzant triangle, a rectangle, or a Hirzebruch trapezoid by a finite number of corner chops.
\end{Theorem}

\begin{proof} Let $\De$ be any Delzant polygon with $d$ edges, let $(v_0, \ldots, v_{d-1})\in\big(\Z^2\big)^d$ be the associated fan of inwards pointing primitive normal vectors, and let $(a_0, \ldots, a_{d-1})\in\Z^d$ be the integers associated to this fan. By Lemma~\ref{lem_fulton} if $d>4$ then $a_i=1$ for some $i\in\{0, \ldots, d-1 \}$ so the fan is the result of a corner chop for some fan with $d-1$ vectors. That is, $\De$ is the result of a corner chop of some Delzant polygon with $d-1$ edges. If $d<5$ then Lemma~\ref{lem_minmodels} lists each possibility. If $d=4$ and $a_0 = a_1 = a_2 = a_3 = 0$ then $\De$ is $\sltz$-equivalent to a rectangle, if $a_0 = 0$, $a_1= k$, $a_2 = 0$, and $a_3 = -k$ for $k\in\Z\setminus\{0\}$ then $\De$ is $\sltz$-equivalent to a Hirzebruch trapezoid, and if $d=3$ with $a_0 = a_1 = a_2 = -1$ then $\De$ is $\sltz$-equivalent to a Delzant triangle.
\end{proof}

\section{Semitoric fans}\label{sec_semitoric}

Now we will apply the method from Section \ref{sec_toric} to classify semitoric fans (Def\/inition \ref{def_stfan}).

The f\/irst step in the classif\/ication is given by a series of lemmas which we will use to manipulate the semitoric fans in a standard form. In Lemma~\ref{lem_cc} we describe the process of corner chopping, in Lemma~\ref{lem_reversecc} we describe the process of reverse corner chopping, in Lemma~\ref{lem_removehidden} we describe the process of removing a~hidden corner, and in Lemma~\ref{lem_fakedelzcommute} we describe the process of commuting a fake and Delzant corner. In each of these lemmas we prove the result of the corresponding process is still a semitoric fan. These four processes are def\/ined in Def\/inition~\ref{def_fantrans}.

\begin{Lemma}\label{lem_cc}
 If $(v_0, \ldots, v_{d-1}) \in\big(\Z^2\big)^d$ is a semitoric fan and $(v_i,v_{i+1})$ is a Delzant corner, then
 \begin{gather*}
 (w_0 = v_0, \ldots, w_i = v_{i}, w_{i+1} = v_i+v_{i+1}, w_{i+2} = v_{i+1},\ldots, w_d = v_{d-1})\in \big(\Z^2\big)^{d+1}
 \end{gather*}
 is a semitoric fan with the pairs $(w_i, w_{i+1})$ and $(w_{i+1},w_{i+2})$ being labeled as Delzant corners.
\end{Lemma}

\begin{proof} It is immediate that $\det(w_i, w_{i+1}) = \det(w_{i+1},w_{i+2})=1$ and the vectors are still in counter-clockwise order because the new vector was inserted between two adjacent vectors.
\end{proof}

\begin{Lemma}\label{lem_reversecc} If $(v_0, \ldots, v_{d-1}) \in\big(\Z^2\big)^d$ is a semitoric fan, $(v_{i-1},v_i)$ and $(v_i, v_{i+1})$ are Delzant corners, and $v_i = v_{i-1} + v_{i+1}$, then
 \begin{gather*}
 (w_0 = v_0, \ldots, w_{i-1} = v_{i-1}, w_i = v_{i+1}, \ldots, w_{d-2} = v_{d-1})\in \big(\Z^2\big)^{d-1}
 \end{gather*}
 is a semitoric fan with the pair $(w_{i-1}, w_i)$ being labeled as a Delzant corner.
\end{Lemma}

\begin{proof} Notice $\det(w_{i-1}, w_i) = \det(v_{i-1}, v_{i} - v_{i-1}) = 1$ so we must only show that the vectors in the new fan are labeled in counter-clockwise order, but this is the case because the vectors in the original fan were and to get the new fan we have only removed a vector.
\end{proof}

\begin{Lemma}\label{lem_removehidden} If $(v_0, \ldots, v_{d-1}) \in\big(\Z^2\big)^d$ is a semitoric fan and $(v_i, v_{i+1})$ is a hidden corner, then
\begin{gather*}
 (w_0=v_0, \ldots, w_i = v_i, w_{i+1}=Tv_{i+1}, w_{i+2}=v_{i+1}, \ldots, w_{d}= v_{d-1})\in\big(\Z^2\big)^{d+1}
\end{gather*}
is a semitoric fan in which $(w_i, w_{i+1})$ is a Delzant corner and $(w_{i+1}, w_{i+2})$ is a fake corner.
\end{Lemma}

\begin{proof} We know that $\det(v_i, T v_{i+1})=1$ because that pair of vectors forms a hidden corner. Notice that \begin{gather*} \det(w_i, w_{i+1}) = \det( v_i, T v_{i+1}) = 1,\end{gather*}
 so that corner is Delzant and
 \begin{gather*}
 \det(w_{i+1}, Tw_{i+2}) = \det(Tv_{i+1,} Tv_{i+1}) =0,
 \end{gather*}
 so if $(w_{i+1}, w_{i+2})$ is on the top boundary then it is a fake corner. It is on the top boundary because $v_{i+1}$ is in the lower half
 plane, $T_\sltz$ sends the lower half plane to the lower half plane, and $(w_{i+1}, w_{i+2})=(Tv_{i+1}, v_{i+1})$. The vectors in the new fan are still in counter-clockwise order because $w_{i+1}$ is in the lower half plane so $w_{i+2} = T^{-1} w_{i+1}$ is counter-clockwise to it and since $w_i$ and $w_{i+1}$ are both in the lower half plane and $\det(w_i, w_{i+1})>0$ we see they are in counter-clockwise order.
\end{proof}

\begin{Lemma}\label{lem_fakedelzcommute} If $(v_0, \ldots, v_{d-1}) \in\big(\Z^2\big)^d$ is a semitoric fan and $(v_i, v_{i+1})$ is a~fake corner and $(v_{i+1}, v_{i+2})$ is a Delzant corner on the top boundary, then
 \begin{gather*}
 (w_0=v_0, \ldots, w_i = v_i, w_{i+1}=Tv_{i+2}, w_{i+2}=v_{i+2}, \ldots, w_{d-1}= v_{d-1})\in\big(\Z^2\big)^{d}
 \end{gather*}
 is a semitoric fan in which $(w_i, w_{i+1})$ is a Delzant corner and $(w_{i+1}, w_{i+2})$ is a fake corner.
\end{Lemma}

\begin{proof}
 We know $\det(v_i, Tv_{i+1}) = 0$ so $v_i = Tv_{i+1}$ since they are both on the top boundary and we also know $\det(v_{i+1}, v_{i+2}) = 1$. Now we can check that
 \begin{gather*}\det(w_i, w_{i+1}) = \det(v_i, Tv_{i+2}) = \det(Tv_{i+1}, Tv_{i+2}) = \det(v_{i+1},v_{i+2}) = 1\end{gather*}
 and
 \begin{gather*}\det(w_{i+1}, Tw_{i+2}) = \det(Tv_{i+2}, Tv_{i+2})=0.\end{gather*}
 Since $(v_{i+1}, v_{i+2})$ is on the top boundary, $v_{i+1}$ and $Tv_{i+1}$ are in the lower half plane, so $(w_{i+1}, w_{i+2}) = (Tv_{i+2}, v_{i+2})$ is in the top boundary which means it is a fake corner, as desired. The argument for why the new vectors are still in counter-clockwise order is similar to the argument in the proof of Lemma~\ref{lem_removehidden}.
\end{proof}

Recall that $\mathcal{G}' = \big\{T^k_\sltz\,|\, k\in\Z\big\}$ is the symmetry group of semitoric fans, given in Def\/ini\-tion~\ref{def_stsymmetry}.

\begin{Lemma}\label{lem_rightangle} Suppose that $(v_0, \ldots, v_{d-1})\in\big(\Z^2\big)^d$ is a semitoric fan. Then after a finite number of corner choppings the fan will be $\mathcal{G}'$-equivalent to one in which two adjacent vectors are $\left(\begin{smallmatrix} 0\\ -1\end{smallmatrix}\right)$ and $\left(\begin{smallmatrix} 1\\ 0\end{smallmatrix}\right)$.
\end{Lemma}

\begin{proof} Let $v_d=v_0$. If \begin{gather*}v_i = \vect{1}{0}\end{gather*} for some $i\in\{0, \ldots, d-1\}$ then notice \begin{gather*}v_{i-1}=\vect{a}{-1}\end{gather*} for some $a\in\Z$. This is because $(v_{i-1}, v_i)$ is not on the top boundary so it must be a~Delzant corner. Then by the action of $T^{-a}_\sltz\in\mathcal{G}'$, which does not change $v_i$, we can attain the required pair of vectors.

Otherwise, renumber so that $v_0$ is in the lower half plane and $v_1$ is in the upper half plane. Then insert the vector $v_0+v_1$ between them. This new vector will have a second component with a strictly smaller magnitude than that of $v_0$ or of $v_1$. Since the magnitude is an integer, repeat this process until the new vector lies on the $x$-axis. Since it is a primitive vector it must be $\left(\begin{smallmatrix} \pm 1\\ 0\end{smallmatrix}\right)$. However since it is the sum of two vectors of opposite sides of the $x$-axis with the one above being counterclockwise about the origin of the one on bottom, it must be~$\left(\begin{smallmatrix} 1\\ 0\end{smallmatrix}\right)$.
\end{proof}

\begin{figure}[t]\centering
 \includegraphics[width=370pt]{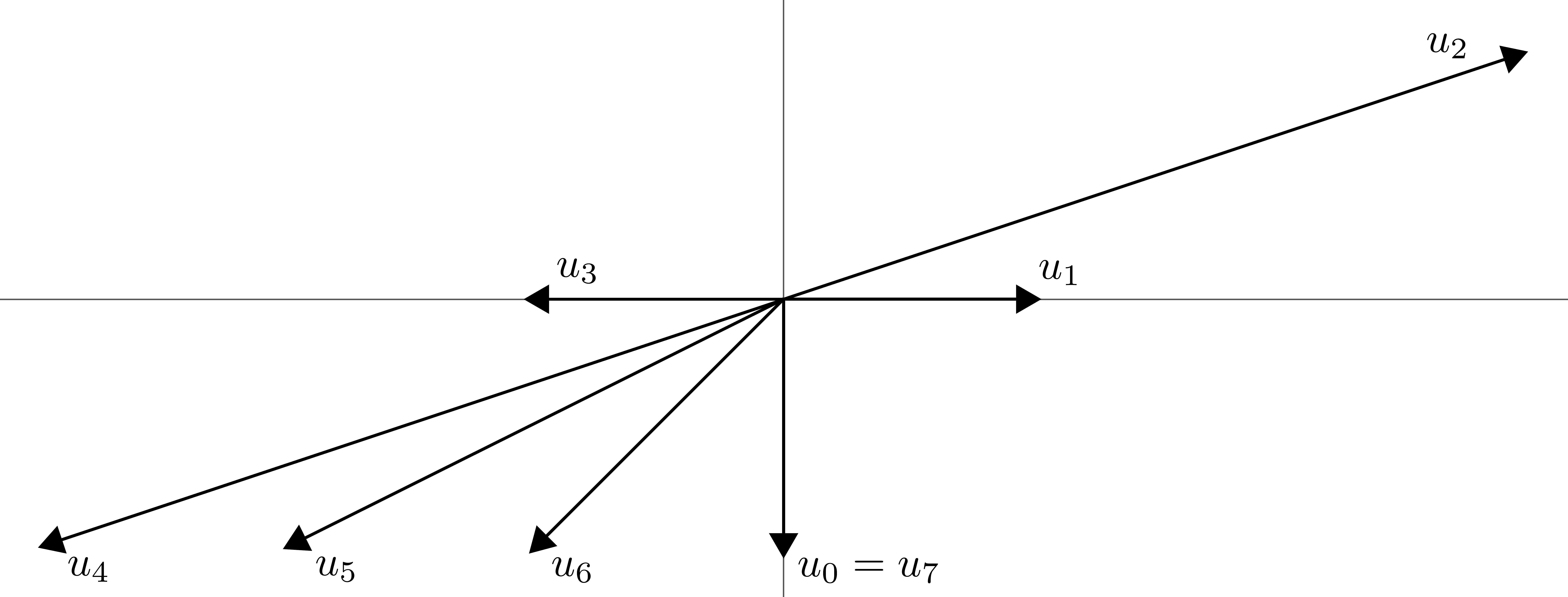}
 \caption{Any semitoric fan with defect $c\in\Z_{\geq0}$ can be transformed into the standard fan of defect $c$. This image has $c=3$.}
 \label{fig_stdfan}
\end{figure}

Now we can put Lemmas \ref{lem_removehidden}, \ref{lem_fakedelzcommute}, and \ref{lem_rightangle} together to produce a standard form for semitoric fans (see Fig.~\ref{fig_stdfan}). This standard form will be important to us in Section \ref{sec_stsyst} because it can be obtained from any semitoric fan of defect~$c$ by only using transformations which are continuous in the space of semitoric polygons.

\begin{Proposition}\label{prop_standardformstfan}
 Let $(v_0, \ldots, v_{d-1})\in\big(\Z^2\big)^d$ be a semitoric fan of defect $c\in\Z_{\geq 0}$.
 \begin{enumerate}[font=\normalfont]\itemsep=0pt
 \item By only corner chopping, reverse corner chopping, removing hidden corners, commuting fake and Delzant corners, and acting by the symmetry group $\mathcal{G}'$ we can obtain a new semitoric fan $(w_0=w_{\ell+c}, \ldots, w_{\ell+c-1})\in\big(\Z^2\big)^{\ell+c}$ with $\ell+c\geq d$ such that
 \begin{itemize}
 \item $w_0 = \left(\begin{smallmatrix} 0\\ -1\end{smallmatrix}\right)$ and $w_1 = \left(\begin{smallmatrix} 1\\ 0\end{smallmatrix}\right)$;
 \item each corner $(w_i, w_{i+1})$ for $i = 0,\ldots, \ell-1$ is Delzant;
 \item each corner $(w_i, w_{i+1})$ for $i= \ell, \ldots, \ell+c-1$ is fake; and
 \item $w_\ell = T^c w_0$ so $\det(w_{\ell-1}, T^c w_0)=1$.
 \end{itemize}
 \item 
 The collection of vectors $(w_\ell, w_1, \ldots, w_{\ell-1})$ is a toric fan;

\item 
The semitoric fan $(w_0, \ldots, w_{c+\ell-1})$ can be obtained via a finite number of corner chops and reverse corner chops from a fan $(u_0, \ldots, u_{c+3})\in\big(\Z^2\big)^{c+4}$ where
 \begin{gather*}
 u_0 = \vect{0}{-1}, \qquad u_1 = \vect{1}{0},\qquad u_2 = \vect{c}{1}, \qquad u_3 = \vect{-1}{0},
 \end{gather*}
 and
 \begin{gather*}
 u_{4+n} = \vect{-c+n}{-1}
 \end{gather*}
 for $n = 0, \ldots, c-1$. In this fan the first four pairs are Delzant corners and the rest are fake corners.
 \end{enumerate}
\end{Proposition}

\begin{proof} The f\/irst part is immediate from Lemmas~\ref{lem_removehidden}, \ref{lem_fakedelzcommute}, and~\ref{lem_rightangle}. By Lemma~\ref{lem_rightangle} we know after a f\/inite number of corner choppings and renumbering it can be arranged that the f\/irst two vectors in the fan are
 \begin{gather*}
 \vect{0}{-1}\qquad \textrm{and} \qquad \vect{1}{0}.
 \end{gather*}
 Then we invoke Lemma \ref{lem_removehidden} to remove all of the hidden corners, and f\/inally use Lemma \ref{lem_fakedelzcommute} commute all of
 the fake corners to be adjacent and arrive at the fan $(w_0, \ldots, w_{\ell+c-1})\in\big(\Z^2\big)^{\ell+c}$ in the statement of the proposition. Notice that $(w_i, w_{i+1})$ being fake for $i= \ell, \ldots, \ell+c-1$ implies that
 \begin{gather*}
 \det(w_\ell, Tw_{\ell+1}) = \cdots = \det(w_{\ell+c-1}, Tw_0) = 0.
 \end{gather*}
 Now, we know both vectors in a fake corner must have negative second component by def\/inition, so this implies that
 \begin{gather*}
 w_\ell = T w_{\ell+1} = T^2 w_{\ell+2} = \cdots = T^c w_0.
 \end{gather*}

 Now $(T^c w_0 = w_\ell, w_1, \ldots, w_{\ell-1})\in\big(\Z^2\big)^\ell$ is a toric fan (we know the vectors are in counter-clockwise order because we started with a semitoric fan) so we have established part~(2) of the theorem. Denote the associated integers for this toric fan by $b_0, \ldots, b_{\ell-1}\in\Z$ so that
 \begin{gather}\label{eqn_propstandardformstfan}
 ST^{b_0}\cdots ST^{b_{\ell-1}} =_G S^4.
 \end{gather}
To prove part~(3) we will use corner chops and reverse corner chops on this toric fan to remove vectors in the set $(w_2, \ldots, w_{\ell-1})$ until the remaining vectors are in the required form. The idea here is that any corner chop or reverse corner chop on the toric fan $(w_\ell, w_1, \ldots, w_{\ell-1})$ which does not remove $w_\ell$ or $w_1$ and does not add a~vector between $w_\ell$ and $w_1$ can also be performed on the semitoric fan $(w_0, \ldots, w_{\ell+c-1})$.

Here we proceed by cases on the values of $\ell$ and $b_0, \ldots, b_{\ell-1}$. We will see in Case~\ref{case_40000} that the semitoric fan is already of the required form and we will show any that other case can be transformed into Case~\ref{case_40000} by corner chops and reverse corner chops on the semitoric fan $(w_0, \ldots, w_{\ell+c-1})$. Here it is important to recall that $w_\ell$ and~$w_1$ are already f\/ixed, so the integers $b_0, \ldots, b_{\ell-1}$ completely determine the toric fan $(w_\ell, w_1, \ldots, w_{\ell-1})$.

\begin{case}\label{case_40000}$\ell=4$ and $b_0=b_1=b_2=b_3=0$.
\end{case}
In this case the toric fan $(w_\ell, w_1, w_2, w_3)$ must be
\begin{gather*}
 \left( \vect{-c}{-1}, \vect{1}{0}, \vect{c}{1}, \vect{-1}{0} \right),
\end{gather*}
which means that the semitoric fan $(w_0, \ldots, w_{\ell+c-1})$ is in the required form for part~(3) of the theorem.

\begin{case}\label{case_4notzero}$\ell = 4$ and $b_0$, $b_1$, $b_2$, $b_3$ not all zero.
\end{case}
In this case, by Lemma~\ref{lem_minmodels}, we must either have $b_0 = b_2 = 0$ and $b_1 = -b_3$, or $b_0 = -b_2$ and $b_1 = b_3 = 0$. First assume $b_0 = b_2 = 0$ and $b_1 = -b_3$. In this case equation~\eqref{eqn_propstandardformstfan} becomes
\begin{gather*}
 S^2T^{b_1}S^2T^{-b_1}=_G S^4.
\end{gather*}
If $b_1>0$ then a corner chop will produce $\big(S^2T^{b_1}STSTST^{-b_1+1}\big)_G$ and a reverse corner chop will produce from this $\big(S^2 T^{b_1-1}S^2T^{-b_1+1}\big)_G$. This can be repeated $b_1$ times to recover Case~\ref{case_40000}. If $b_1<0$ then a corner chop can produce
$\big(S^2T^{b_1+1}STSTST^{-b_1}\big)_G$ and from this a reverse corner chop can produce $\big(S^2T^{b_1+1}S^2T^{-b_1-1}\big)_G$. Again, this can be repeated $b_1$ times to recover Case~\ref{case_40000}. Observe that the above corner chops and reverse corner chops do not af\/fect the f\/irst and last element of the toric fan; thus the above operations can be realized as corner chops and reverse corner chops on the semitoric fan $(w_0, \ldots, w_{\ell+c-1})$ as required.

The case of $b_0 = -b_2$ and $b_1 = b_3 = 0$ is similar.

\begin{case}\label{case_ellless4} $\ell<4$.
\end{case}
 By Lemma~\ref{lem_minmodels} we must have $\ell=3$. Perform a corner chop between the vectors $w_1$ and $w_2$ to add the vector $w_1+w_2$ and reduce to either Case~\ref{case_40000} or Case~\ref{case_4notzero}.

\begin{case}\label{case_ellmore4} $\ell>4$.
\end{case}
In this case we will show that the semitoric fan can always be transformed into another semitoric fan with one less length, so repeating this
process will eventually yield $\ell=4$, thereby reducing to Case~\ref{case_40000} or Case~\ref{case_4notzero}. There are two subcases in this case.

\begin{subcase}\label{case_revcornerchop} $b_i=1$ for some $i$ such that $0<i<\ell-1$.
\end{subcase}
 Here a reverse corner chop may be performed on $w_{i+1}$, removing it from the semitoric fan and reducing the length by~1.

\begin{subcase}\label{case_lastcase} $b_i \neq 1$ for all $i=1, \ldots, \ell-2$.
\end{subcase}
 Any substitution using $STS =_G T^{-1}ST^{-1}$ in equation~\eqref{eqn_propstandardformstfan} can be realized as a corner chop or reverse corner chop on the toric fan $(w_\ell, w_1, \ldots, w_{\ell-1})$ and no such substitution corresponds to the corner chop on $(w_{\ell-1}, w_\ell)$ so these transformations can also be realized a corner chops and reverse corner chops on the semitoric fan $(w_0, \ldots, w_{\ell+c-1})$. In particular, if $\big(ST^{b_0}\cdots ST^{b_{\ell-1}}\big)_G$ contains the subword $\big(T^{a+1}SST^b\big)_G$, we can perform a corner chop to obtain $\big(T^{a+1}STSTST^{b+1}\big)_G$, and then a reverse corner chop to reduce to $\big(T^aSST^{b+1}\big)_G$. Note that this can be done even if the $T^a_G$ was in the $\big(ST^{b_0}\big)_G$ term and also notice that a similar argument can be used to transform $\big(T^{a-1}SST^b\big)_G$ into $\big(T^{a}SST^{b-1}\big)_G$ if $a<0$. This means that every occurrence of $S^2_G$ may be commuted to the front of the word, and thus we may assume that there is an integer $m\geq0$ for which $b_i = 0$ if and only if $i<m$. If $m\geq 4$ then
 \begin{gather*}
 \big(ST^0\big)^{m-4}ST^{b_m}\cdots ST^{b_{\ell-1}}=_G I,
 \end{gather*}
 so
 \begin{gather*}
 w_G\big(\big(ST^0\big)^{m-4}ST^{b_m}\cdots ST^{b_{\ell-1}}\big) = w_G(I) = 0,
 \end{gather*}
 but the integers $(0,\ldots,0,b_m, b_{m+1}, \ldots, b_{\ell-1})$ have associated vectors which travel some integer $k>0$ times around the origin, since each vector is counterclockwise from the previous one, which means that
 \begin{gather*}
 w_G\big(\big(\big(ST^0\big)^{m-4}ST^{b_m}\cdots ST^{b_{\ell-1}}\big)_G\big) = 12k
 \end{gather*}
 forming a contradiction. Thus, $m<4$ and this also implies $\ell-m>1$ because $\ell>4$.

Let $e=0$ if $m$ is even and $e=1$ if $m$ is odd, then equation~\eqref{eqn_propstandardformstfan} implies the following equality in $\psl$:
 \begin{gather}\label{eqn_propststandard2}
 S^{e} ST^{b_m}\cdots ST^{b_{\ell-1}}=_{\psl} I.
 \end{gather}
 We proceed by cases on the value of $\ell-m+e$.

 We already know that $\ell-m+e\geq \ell-m>1$, and also that $\ell-m+e = 2$ is impossible by applying Lemma~\ref{lem_psltz} to equation~\eqref{eqn_propststandard2}. If $e=1$ then $\ell-m+e=3$ is impossible by Lemma~\ref{lem_psltz} and if $e=0$ then $\ell-m+e=3$ is impossible by acting with $w_G$ on both sides of equation~\eqref{eqn_propstandardformstfan}. By Lemma~\ref{lem_psltz} if $e=0$ and $\ell-m+e=3$ then equation~\eqref{eqn_propstandardformstfan} must be
 \begin{gather*}
 S^mST^{-1}ST^{-1}ST^{-1} =_G S^4,
 \end{gather*}
 which implies
 \begin{gather*}
 w_G\big(\big(S^mST^{-1}ST^{-1}ST^{-1}\big)_G\big) = w_G\big( S^4_G\big),
 \end{gather*}
 which can only hold if $m=0$ in which case $\ell-m+e=3$ implies that $\ell =3$, which is false. If $\ell-m+e=4$ then either $b_{\ell-1}=0$ or $b_{\ell-2}=0$, both of which are false by assumption.

 The only remaining case is $\ell-m+e >4$, and in this case Lemma~\ref{lem_psltz} implies that there exists some $i$ with $m< i <\ell-1$ and $b_{i}=\{-1, 0, 1\}$. Thus, $b_i = \pm 1$. If $b_i=1$, then we may perform a reverse corner chop on $w_{i+1}$, as in Case~\ref{case_revcornerchop},
 to reduce $\ell$. Otherwise, $b_i = -1$ and we use the relation $ST^{-1}S=_G SSTST$. Since
 \begin{gather*}
 ST^{-1}S =_G S\big(T^{-1}ST^{-1}\big)T =_G S(STS)T=_G SSTST
 \end{gather*}
 this relation can actually be realized by a corner chop, which we know corresponds to a legal move at the level of fans. Here we have introduced an $S^2_G$ term, and we next return to the beginning of Case~\ref{case_lastcase}. That is, the $S^2_G$ term is moved to the front of the word and now the f\/irst~$m'$ powers of~$T_G$ are zero, where $m'>m$, and the new toric fan has length $\ell' >\ell$. If the new word contains $(STS)_G$ then a reverse corner chop may be performed, as in Case~\ref{case_revcornerchop}, before returning to the beginning of Case~\ref{case_lastcase} (with a~larger $m$ and $\ell$) and otherwise the word contains $\big(ST^{-1}S\big)_G$, which we can again remove. This process must terminate by f\/inding an occurrence of $(STS)_G$ and performing a reverse corner chop after at most 3 iterations because at each iteration the integer $m$ increases and we have already argued that $m\geq4$ is impossible. Once $m=3$ then no iteration of this process can include $\big(ST^{-1}S\big)_G$ because in such a~case~$m$ can be increased by one. This is important because on those iterations during which the replacement $\big(ST^{-1}S\big)_G$ to $\big(S^2TST\big)_G$ is used the total number $\ell$ of vectors in the toric fan increases by one, but we have shown this can happen at most three times.

 Thus, the algorithm described in Case~\ref{case_lastcase} must eventually allow for enough reverse corner chops to be performed to reducing the number of vectors in the toric fan to four. This completes Case~\ref{case_lastcase}.

 Therefore, in any case with $\ell\neq 4$ the semitoric fan can be transformed into one with $\ell=4$ by Cases~\ref{case_ellless4} and~\ref{case_ellmore4}, and any fan with $\ell=4$ is either in the correct form already (Case~\ref{case_40000}) or can be transformed into the correct form (Case~\ref{case_4notzero}).
\end{proof}

Theorem \ref{thm_stpoly} is immediate from Proposition~\ref{prop_standardformstfan}.

\begin{Remark} Notice that Theorem \ref{thm_toricpoly} is dif\/ferent from Theorem~\ref{thm_stpoly} because in Theo\-rem~\ref{thm_toricpoly} the minimal models of the Delzant polygons may be achieved through only corner chops. In Theorem~\ref{thm_stpoly} we use instead a variety of transformations and in Section~\ref{sec_stsyst} we show all of those transformations are continuous with respect to the topology on the moduli space of semitoric systems (described in Section~\ref{sec_stmetric}).
\end{Remark}

\newcommand{\rxyz}{\mathbb{R}[[X,Y]]_0}
\newcommand{\bn}{\{b_n\}_{n=0}^\infty}
\newcommand{\dts}{d_{\rxyz}}

\newcommand{\poly}{{\rm Polyg}\big(\R^2\big)}
\newcommand{\lwpoly}{{\rm LW}{\rm Polyg}\big(\R^2\big)}

\newcommand{\stpoly}{\mathcal{D}{\rm Polyg}\big(\R^2\big)}
\newcommand{\stpolymfk}{\mathcal{D}{\rm Polyg}_{\mf, \vec{k}}\big(\R^2\big)}

\newcommand{\vertr}{{\rm Vert}\big(\R^2\big)}
\newcommand{\lwp}{\big(\De, (\ell_{\la_j}, \ep_j, k_j)_{j=1}^{\mf}\big)}
\newcommand{\pstpoly}{\big(\De, (\ell_{\la_j}, +1, k_j)_{j=1}^{\mf}\big)}
\newcommand{\pstpolyprime}{\big(\De', (\ell_{\la'_j}, +1, k'_j)_{j=1}^{\mf}\big)}
\newcommand{\symdiff}{\ast}
\newcommand{\distpoly}{d_{\mathcal{P}}^{\nu}}
\newcommand{\topbound}{\partial^{{\rm top}}\De}

\newcommand{\mf}{{m_f}}

\newcommand{\m}{\mathcal{M}}
\newcommand{\distm}{d^{\nu, \bn}}
\newcommand{\mmf}{\mathcal{M}_{\mf}}
\newcommand{\mmfk}{\mathcal{M}_{\mf, \vec{k}}}
\newcommand{\tmfk}{\mathcal{T}_{\mf, \vec{k}}}
\newcommand{\distmmfk}{d^{\nu, \bn}_{\mf, \vec{k}}}
\newcommand{\melement}{\big([\De_w], (h_j)_{j=1}^{\mf}, ((S_j)^\infty)_{j=1}^{\mf}\big)}
\newcommand{\melementprime}{\big([\De'_w], (h'_j)_{j=1}^{\mf}, ((S'_j)^\infty)_{j=1}^{\mf}\big)}
\newcommand{\melementallprime}{\big([\De'_w], (h'_j)_{j=1}^{\mf'}, ((S'_j)^\infty)_{j=1}^{\mf'}\big)}

\newcommand{\distt}{\mathcal{D}^{\nu, \bn}}

\section{Application to symplectic geometry} \label{sec_stsyst}

The results of Section \ref{sec_semitoric} have an interpretation in the symplectic geometry of toric manifolds and semitoric integrable systems. In Section \ref{sec_invariantsofst} we will review the classif\/ication theorem of Pelayo--V\~{u} Ng\d{o}c~\cite{PeVNconstruct2011}, in Section~\ref{sec_stmetric} we def\/ine the metric on the moduli space of semitoric systems given by the second author in~\cite{PaSTMetric2015}, and in Section~\ref{sec_stconnect} we prove the connectivity result for the moduli space of semitoric systems (with respect to the aforementioned metric) using Theorem~\ref{thm_stpoly}.

\begin{Definition}[\cite{PeVNsemitoricinvt2009}]\label{def_stsystem} A \emph{semitoric integrable system} is a 4-dimensional connected symplectic
 manifold $(M,\om)$ with an integrable Hamiltonian system $F=(J,H)\colon M\to\R^2$ such that $J$ is a proper momentum map for an ef\/fective Hamiltonian $S^1$-action on~$(M,\om)$ and $F$ has only non-degenerate singularities which have no real-hyperbolic blocks. Such a system is said to be a~\emph{simple semitoric integrable system} if there is at most one focus-focus critical point in $J^{-1}(x)$ for any $x\in\R$.

 An \emph{isomorphism of semitoric systems} is a symplectomorphism $\phi\colon (M_1,\om_1) \to (M_2,\om_2)$, where $(M_1, \om_1, F_1=(J_1,H_1))$ and
 $(M_2, \om_2, F_2=(J_2,H_2))$ are semitoric systems, such that $\phi^* (J_2, H_2) = (J_1, f(J_1,H_1))$ where $f\colon F_1(M_1)\to\R$ is a smooth function such that $\deriv{f}{H_1}$ is eve\-ry\-where positive. We denote the moduli space of simple semitoric systems modulo semitoric isomorphism by $\mathcal{T}$.
\end{Definition}

\subsection{Invariants of semitoric systems}\label{sec_invariantsofst}

The momentum map image of a compact toric integrable system, which is a polytope, is suf\/f\/icient to classify such systems up to a suitable
notion of isomorphism (cf.~\cite{De1988}). Semitoric systems are classif\/ied in terms of a list of invariants~\cite{PeVNconstruct2011}. Roughly speaking, the complete invariant is a polygon together with a set of interior points each labeled with extra information (encoding singularities of so called \emph{focus-focus type}, which semitoric systems may possess, but toric systems do not) modulo an equivalence relation. Even without the extra information semitoric polygons are more complicated than toric polygons because they have fake corners and hidden corners (as def\/ined in Def\/inition~\ref{def_corners}). In this section we will def\/ine each of the invariants of semitoric systems and give the necessary
def\/initions to state the Pelayo--V\~{u} Ng\d{o}c classif\/ication theorem. Readers interested in further details may
consult~\cite{PeVNsemitoricinvt2009, PeVNconstruct2011,PeVNsymplthy2011, PeVNfirst2012, VN2003, VN2007}.

\subsubsection{The number of focus-focus points invariant}\label{sec_numberffpoints}

While a toric integrable system can only have transversally-elliptic and elliptic-elliptic singularities a semitoric integrable system can also have focus-focus singularities (for a def\/inition of these types of singularities see for instance~\cite{PeVNconstruct2011}). In~\cite[Theorem~1]{VN2007} V\~{u} Ng\d{o}c proves that any semitoric system has at most f\/initely many focus-focus singular points. The f\/irst invariant is a~nonnegative integer $\mf$ known as the \emph{number of focus-focus singular points invariant}.

\subsubsection{The Taylor series invariant}

The semi-global\footnote{I.e.,~in a neighborhood of the f\/iber over the critical point.} structure of a focus-focus singular point is determined by a~formal power series in two variables up to the suitable notion of isomorphism (see~\cite{VNSnotes,VN2003}).

\begin{Definition} Let $\mathbb{R}[[X,Y]]$ denote the algebra of real formal power series in two variables and let $\rxyz \subset \mathbb{R}[[X,Y]]$ be the subspace of series $\sum\limits_{i,j\geq0}\si_{i,j}X^iY^j$ which have $\si_{0,0}=0$ and $\si_{0,1}\in[0,2\pi)$.
\end{Definition}

The \emph{Taylor series invariant} is one element of $\rxyz$ for each of the $\mf$ focus-focus points. Thanks to~\cite{VNSnotes} we understand why the isomorphism need not be taken into account when there is a global $S^1$-action, that is, as in the case of semitoric systems~\cite{PeVNsemitoricinvt2009}, provided one assumes everywhere that the Eliasson isomorphisms preserve the global $S^1$-action and the $\R^2$-orientation, in which case the Taylor series is unique (for the general case uniqueness is up to a~$(\Z_2\times\Z_2)$-action, see~\cite{VNSnotes}).

\subsubsection{The polygon invariant and the twisting index invariant}\label{sec_affinetwist}

The polygon invariant is the invariant which is directly analogous to the Delzant polygon in the toric case and the twisting index is an integer label on each focus-focus f\/iber which encodes how the semiglobal (in a neighborhood of a f\/iber) model of the focus-focus point relates to the global toric momentum map used to generate the polygon. The dif\/ferent choices of toric momentum map produce a family of polygons and thus the component of the twisting index associated to a~single focus-focus f\/iber may change depending on this choice, but if the system has more than one focus-focus point then the dif\/ference of twisting index components of focus-focus f\/ibers does not depend on the choice of polygon. This is because dif\/ferent choices of polygon can only change the sequence of integers which is the twisting index by the addition of a~common integer, and for this reason we are able to also def\/ine the twisting index class (f\/irst introduced in~\cite[Def\/inition 3.7]{PaSTMetric2015}, and def\/ined in the present paper in Def\/inition~\ref{def_twistingindexclass}). The polygon and twisting index invariants are described together because the choice of common integer added to the twisting index is related to the polygon.

In general, the momentum map image of a semitoric system does not have to be a polygon and does not even have to be convex, but in~\cite{VN2007} the author was able to recover a family of convex polygons which take the place of the Delzant polygon in the semitoric case.

\begin{Definition} Let $\poly$ denote the set of all convex polygons (as in Def\/inition~\ref{def_poly}) which are rational (i.e., each edge is directed along a vector with integer coef\/f\/icients), and observe that all elements of $\poly$ are automatically simple (two edges meet at each vertex).
\end{Definition}

For brevity, for the duration of this paper by ``polygon'' we will always mean ``rational, convex polygon''. Notice that a polygon may be noncompact. For any $\la\in\R$ we will use the notation
 \begin{gather*} \ell_\la = \big\{ (x,y)\in\R^2 \,|\, x = \la\big\}\qquad \text{and} \qquad \vertr = \{ \ell_\la \,|\, \la \in \R\}.\end{gather*}

 \begin{Definition} A \emph{labeled weighted polygon of defect\footnote{In this context this is sometimes called the \emph{complexity},
 but seeing as complexity sometimes refers to half the dimension of the manifold minus the dimension of the acting
 torus (so all semitoric systems would be complexity~1) we have chosen to use the word defect.} $\mf\in\Z_{\geq0}$}
 is def\/ined to be
 \begin{gather*} \lwp \in \poly \times \big(\vertr\times\{-1,+1\}\times\Z\big)^{m_f} \end{gather*}
 with
 \begin{gather*} \inf_{s\in\De}\pi_1(s)< \la_1 < \cdots < \la_{\mf} < \sup_{s\in\De} \pi_1 (s),\end{gather*}
 where $\pi_1\colon \R^2 \to \R$ is projection onto the f\/irst coordinate. We denote the space of labeled weighted polygons of any defect by $\lwpoly$.
 \end{Definition}

The labels $k_j$, $j=1, \ldots, \mf$, will become the twisting index invariant after the appropriate quotient. The polygon invariant will be
the orbit of a specif\/ic type of element of $\lwpoly$. Recall\footnote{The authors of~\cite{PeVNsemitoricinvt2009,PeVNconstruct2011} denote by $T$ the transpose of this matrix. Since in this paper we have discovered the connection to $\sltz$ we have chosen to instead use the notation standard to $\sltz$-presentations.}
\begin{gather*}
 T_\sltz = \matr{1}{1}{0}{1}\in\sltz \qquad \text{so} \qquad \big(T^t_\sltz\big)^k = \matr{1}{0}{k}{1}
\end{gather*}
for $k\in\Z$. Recall by Notation~\ref{notationST} that $Tv$ always means $T_\sltz v$ for any $v\in\Z^2$.

 \begin{Definition}\label{def_corners}\quad
 \begin{enumerate}\itemsep=0pt
 \item For $\De\in\poly$ a point $(x_0,y_0)\in\R^2$ is said to be in the \emph{top boundary of $\De$} if $y_0 = \sup \{y\in\R \,|\, (x_0, y) \in \De\}$. We denote this by $(x_0,y_0)\in\topbound$.
 \item Let $\De\in\poly$. A \emph{corner} of $\De$ is a point $p\in\partial\De$ such that the edges meeting at $p$ are not co-linear. For a corner $p\in\De$ let $u,v\in\R^2$ be primitive inwards pointing normal vectors to the edges of $\De$ adjacent to $p$ in the order of positive orientation. We say that $p$ satisf\/ies
 \begin{enumerate}\itemsep=0pt
 \item the \emph{Delzant condition} if $\det(u,v) = 1$;
 \item the \emph{hidden condition} if $p\in\topbound$ and $\det(u,Tv) = 1$; and
 \item the \emph{fake condition} if $p\in\topbound$ and $\det(u, Tv) = 0$.
 \end{enumerate}
 \item We say $\De\in\poly$ has \emph{everywhere finite height} if the intersection of $\De$ with any vertical line is either compact or empty.
 \end{enumerate}
 \end{Definition}

 \begin{Definition}\label{def_semitoricprimitive} An element $\lwp\in\lwpoly$ is a \emph{primitive semitoric polygon} if
 \begin{enumerate}\itemsep=0pt
 \item[1)] $\De$ has everywhere f\/inite height;
 \item[2)] $\ep_j = +1$ for $j=1, \ldots, \mf$;
 \item[3)] any point in $\topbound \cap \ell_{\la_j}$ for some $j\in\{1, \ldots, \mf\}$ is a corner which satisf\/ies either the hidden or fake condition; and
 \item[4)] all other corners satisfy the Delzant condition.
 \end{enumerate}
 \end{Definition}

 As is done in \cite{PeVNconstruct2011} we say that those corners of a primitive semitoric polygon which are in the top boundary and also in some line $\ell_{\la_j}$ are either a \emph{hidden corner} or \emph{fake corner} corresponding to which condition they satisfy. Also, all other corners of a primitive semitoric polygon must by def\/inition satisfy the Delzant condition and thus are referred to as \emph{Delzant corners}. Notice here that the corners of a primitive semitoric polygon do not need to be labeled by their type (hidden, fake, or Delzant) because the type of a corner can by deduced by examining the primitive semitoric polygon (compare to semitoric fans, see Def\/inition \ref{def_stfan} and the following discussion).

 In the case that $\lwp$ is a primitive semitoric polygon and $\De$ is compact then there exists an associated semitoric fan.

 \subsubsection[The action of $G_{\mf}\times\mathcal{G}$]{The action of $\boldsymbol{G_{\mf}\times\mathcal{G}}$}

Now we will def\/ine a group and the way that its elements act on a labeled weighted polygon of defect $\mf$, for some f\/ixed $\mf$.
For any $\ell\in\vertr$ and $k\in\Z$ f\/ix an origin in $\ell$ and let $t_{\ell}^k\colon \R^2\to\R^2$ act as the identity on the half-space to the left of $\ell$ and as $(T^t)^k$, with respect to the origin in $\ell$, on the half-space to the right of $\ell$. For $\vec{u} = (u_1, \ldots, u_{\mf})\in\{-1,0,1\}^{\mf}$ and $\vec{\la} = (\la_1, \ldots, \la_{\mf})\in \R^{\mf}$ let \begin{gather*} t_{\vec{\la}}^{\vec{u}} = t^{u_1}_{\ell_{\la_1}} \circ \cdots \circ t^{u_{\mf}}_{\ell_{\la_{\mf}}}.\end{gather*}

 \begin{Definition}\label{def_actionofGGmf}
 For any nonnegative $\mf\in\Z$ let $G_{\mf} = \{-1,1\}^{\mf}$ and let $\mathcal{G} = \big\{ \big(T^t_\sltz\big)^k \,|\, k\in\Z\big\}$. We def\/ine the action of $((\ep'_j)_{j=1}^{\mf}, (T^t_\sltz)^k)\in G_{\mf}\times\mathcal{G}$ on an element of $\lwpoly$ by \begin{gather*} \big((\ep'_j)_{j=1}^{\mf}, \big(T^t_\sltz\big)^k\big) \cdot \lwp = \big(t^{\vec{u}}_{\vec{\la}}\big(\big(T^t_\sltz\big)^k \De\big), (\ell_{\la_j}, \ep'_j\ep_j, k+k_j)_{j=1}^{\mf}\big),\end{gather*} where $\vec{\la} = (\la_1, \ldots, \la_{\mf})$ and $\vec{u} = \big(\frac{\ep_j - \ep_j \ep'_j}{2}\big)_{j=1}^{\mf}$.
 \end{Definition}

 \begin{Definition}\label{def_stpoly}
 A \emph{semitoric polygon} is the orbit under $G_{\mf}\times\mathcal{G}$ of
 a primitive semitoric polygon. That is, given a primitive semitoric polygon
 $\De_w = \pstpoly\in\lwpoly$ the associated semitoric polygon is the subset
 of $\lwpoly$ given by
\begin{gather*}[\De_w]=\big\{ \big(t_{\vec{\la}}^{\vec{u}}\big(\big(T^t_\sltz\big)^k(\De)\big), (\ell_{\la_j}, 1-2u_j, k_j+k)_{j=1}^{\mf}\big) \,|\, \vec{u}\in\{0,1\}^{\mf}, k\in\Z\big\}.
\end{gather*}
 The collection of semitoric polygons is denoted $\stpoly$.
 \end{Definition}

 In general the action of $G_{\mf}\times \mathcal{G}$ may not preserve the convexity of the polygons but it is shown in \cite[Lemma 4.2]{PeVNconstruct2011} that $[\De_w]\subset\lwpoly$ for any primitive semitoric polygon $\De_w$.

 \begin{Remark}
 A semitoric polygon is a family of polygons which is determined by the choice of a single primitive semitoric polygon, though this choice is not unique since inf\/initely many primitive semitoric polygons generate any given semitoric polygon. For instance, $[\De, (\ell_0,+1,0)]$ and $[\De',(\ell_0,+1,1)]$ are the same semitoric polygon where $\De$ is the convex hull of $(-1,0)$,$(0,1)$,$(1,0)$, and $(2,1)$, and $\De'$ is the convex hull of $(-1,-1)$, $(0,1)$, $(1,1)$, and $(2,3)$, since $T^t(\De)=\De'$.
 \end{Remark}

\begin{Definition}\label{def_primitivetwistingindex}
 Let $\De_w = (\De,(\ell_{\la_j}, \ep_j, k_j)_{j=1}^\mf)$ be a primitive semitoric polygon. Then the \emph{twisting index of $\De_w$} is the vector $\vec{k} = (k_1, \ldots, k_{\mf})\in\Z^\mf$.
\end{Definition}

The twisting index of a semitoric system is the coherent assignment of a twisting index to each element of the associated semitoric polygon, which is to say that the twisting index invariant is the set
\begin{gather*}
 \big[\big(\De,(\ell_{\la_j}, \ep_j, k_j)_{j=1}^\mf\big)\big]\in\lwpoly/(G_{\mf}\times \mathcal{G}).
\end{gather*}
\begin{Remark}
In~\cite{PeVNsemitoricinvt2009} the authors def\/ine the polygon invariant as this set without the integer labeling (without the $k_j$) and the twisting index invariant as this set including the labeling. In this paper we def\/ine these together because this allows us to give a more transparent formulation of Theorem~\ref{thm_B}.
\end{Remark}

\begin{Definition}\label{def_twistingindexclasspoly}
 Two semitoric polygons $[\De_w]$ and $[\De_w']$ are in the same \emph{twisting index class} if and only if the primitive semitoric polygons
 \begin{gather*}\De_w = \big(\De,(\ell_{\la_j}, \ep_j, k_j)_{j=1}^\mf\big)\qquad \textrm{and} \qquad {\De_w}' = \big(\De',(\ell_{{\la_j}'}, {\ep_j}', {k_j}')_{j=1}^{\mf'}\big),\end{gather*}
 representing them can be chosen such that $\mf = \mf'$ and $k_j = k_j '$ for $j=1, \ldots, \mf$.
\end{Definition}

\subsubsection{The volume invariant}

Suppose that $\De_w$ is a primitive semitoric polygon. For each $j=1,\ldots, \mf$ we def\/ine a real number $h_j \in (0, \textrm{length}(\pi_2 (\De \cap \ell_{\la_j})))$ by measuring the Liouville volume of a specif\/ic subset of~$M$ related to each focus-focus point. Specif\/ically, if $p_j\in M$ is a focus-focus point then we def\/ine
\begin{gather*}
 h_j = \operatorname{vol}\{x\in M\,|\, J(x) = J(p_j)\textrm{ and } H(x)<H(p_j)\},
\end{gather*}
where $\operatorname{vol}$ is the Liouville volume on $M$. The real numbers $h_1,\ldots, h_{\mf}$ are independent of the choice of primitive semitoric polygon representing the semitoric polygon associated to a given semitoric system. The details are in~\cite{PeVNsemitoricinvt2009, PeVNconstruct2011}.

\subsubsection{The complete invariant}
\begin{figure}[t]\centering
 \includegraphics[width = 350pt]{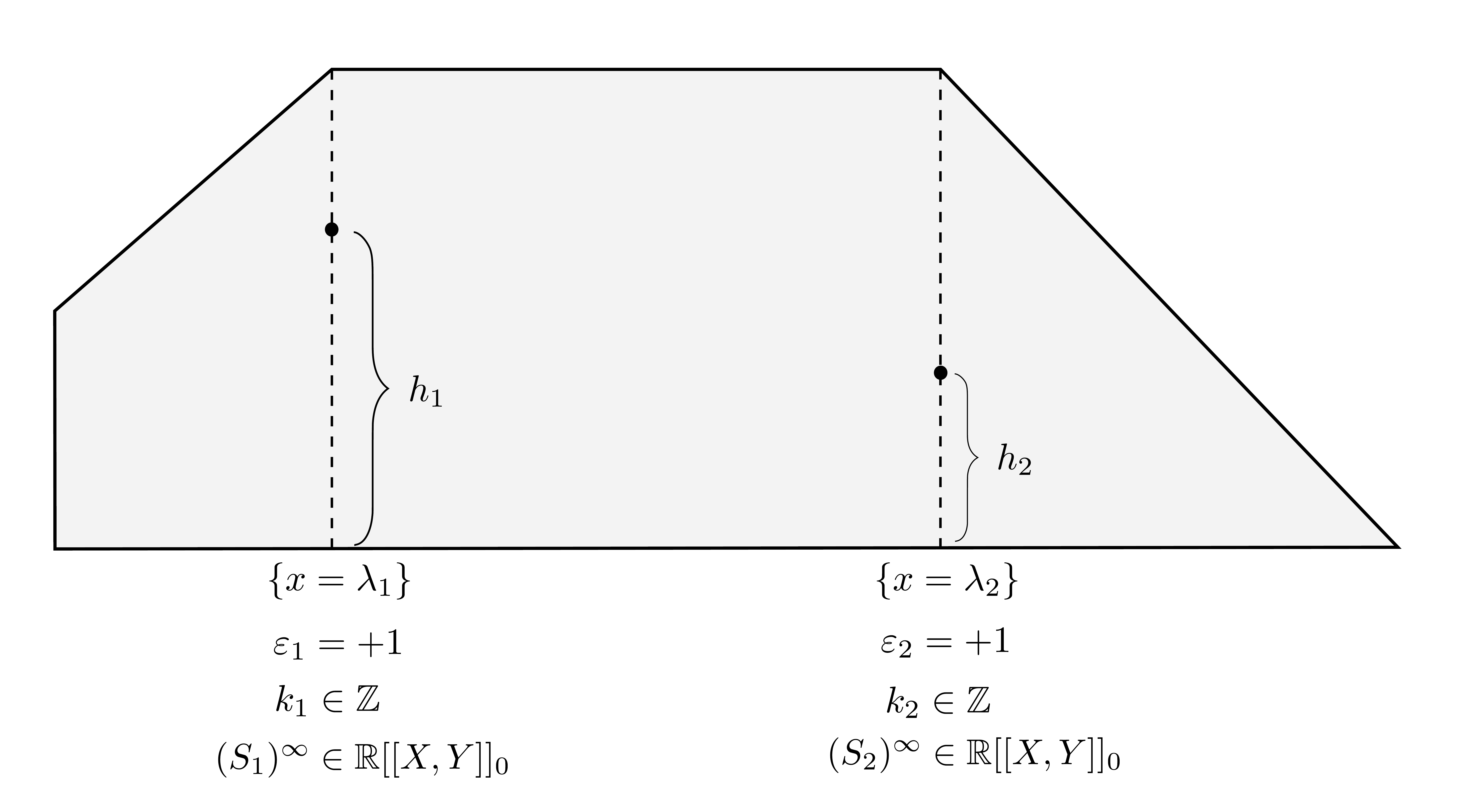}
 \caption{The complete invariant of semitoric systems is the orbit under $G_{\mf}\times\mathcal{G}$ of this object.} \label{fig_completeinvarient}
 \end{figure}

The collection of these f\/ive invariants forms the complete invariant of a semitoric system, which can be represented by a~polygon $\De$ (satisfying the conditions of Def\/inition~\ref{def_semitoricprimitive}) with $\mf$ vertical lines $\ell_{\la_1}, \ldots, \ell_{\la_{\mf}}$ and for each $j=1, \ldots, \mf$ there is a distinguished point on $\ell_{\la_j}$ a~distance~$h_j$ from the bottom of~$\De$. Each of these points is also labeled with an integer $k_j$ and a Taylor series~$(S_j)^\infty$. This is shown in Fig.~\ref{fig_completeinvarient}. It is important to note that in this rough description we have omitted the action of group $G_{\mf}\times\mathcal{G}$ which is a key part to understanding the classif\/ication.

\subsubsection{The classif\/ication theorem}

 \begin{Definition}[\cite{PeVNconstruct2011}]\label{def_listofingre}
 A \emph{semitoric list of ingredients} is
 \begin{enumerate}\itemsep=0pt
 \item[1)] 
 a nonnegative integer $\mf$;
 \item[2)] 
 a semitoric polygon $[\De_w]=[\pstpoly]$ of defect $\mf$;
 \item[3)] 
 a collection of $\mf$ real numbers $h_1, \ldots, h_{\mf} \in \R$ such that $0<h_j<\text{length}(\pi_2(\De\cap\ell_{\la_j}))$ for each $j=1,\ldots,\mf$; and
 \item[4)] a collection of $\mf$ Taylor series $(S_1)^\infty, \ldots, (S_{\mf})^\infty\in\rxyz$.
 \end{enumerate}
 Let $\m$ denote the collection of all semitoric lists of ingredients and let $\mmf$ be lists of ingredients with ingredient~(1) equal to the nonnegative integer $\mf$.
 \end{Definition}

We use the following result to study the moduli space of semitoric systems by instead studying the set of semitoric lists of ingredients.

 \begin{Theorem}[\cite{PeVNconstruct2011}]\label{thm_class}
 There exists a bijection between the set of simple semitoric integrable systems modulo semitoric isomorphism and $\m$, the set of semitoric lists of ingredients. In particular,
 \begin{align*}
 \mathcal{T} &\cong \m,\\
 [(M,\om,(J,H))] &\leftrightarrow \big([\De_w], (h_j)_{j=1}^{\mf}, ((S_j)^\infty)_{j=1}^{\mf}\big),
 \end{align*}
 where the invariants $\mf$, $[\De_w]$, $h_j$, and $(S_j)^\infty$ are as defined above, and the bijection is as given in~{\rm \cite{PeVNconstruct2011}}.
 \end{Theorem}

Given a semitoric system, the invariants can be obtained as outlined above, and given semitoric list of ingredients a semitoric system can be constructed by a gluing procedure. Explicitly constructing the invariants from the system and explicitly constructing the system from a given set of invariants is described in detail in~\cite{PeVNsemitoricinvt2009, PeVNconstruct2011}.

\begin{Definition}\label{def_twistingindexclass} Two semitoric systems $(M,\om,F)$ and $(M,\om,F')$ are said to be in the same \emph{twisting index class} if the semitoric polygon associated to $(M,\om,F)$ and the semitoric polygon associated to $(M',\om',F')$ are in the same twisting index class, as in Def\/inition~\ref{def_twistingindexclasspoly}.
\end{Definition}

\subsection{Metric and topology on the moduli space}\label{sec_stmetric}

In~\cite{PaSTMetric2015} the second author def\/ines a metric space structure for $\mathcal{T}$ using Theorem \ref{thm_class}. This metric is constructed by combining metrics on each ingredient. We reproduce this construction brief\/ly in this section.

\subsubsection{The metric on Taylor series}\label{sec_metricts}

Recall that the Taylor series in $\rxyz$ have the coef\/f\/icient to $X^0Y^1$ in the interval $[0,2\pi)$. From the construction in~\cite{VN2003} we
can see that when def\/ining the topology the endpoints of the closure of this interval should be identif\/ied. We def\/ine a topology on $\rxyz$ such that a~sequence of Taylor series converges if and only if each term converges in the appropriate space. This is the topology induced by any of the following metrics.

 \begin{Definition} A sequence $\bn$ with $b_n \in (0, \infty)$ for each $n\in\Z_{\geq0}$ is called \emph{linearly summable} if $\sum\limits_{n=0}^\infty\! n b_n <\infty$. Let $\bn$ be such a sequence and def\/ine \smash{$\dts\!\colon (\rxyz)^2 \!\to\! \R$} by
 \begin{gather*}
 \dts((S)^\infty,(S')^\infty) = \sum_{i,j\geq 0,(i,j)\neq(0,1)}^\infty \min\{\abs{\si_{i,j}-\si'_{i,j}},b_{i+j}\}\\
 \hphantom{\dts((S)^\infty,(S')^\infty) =}{} + \min\{\abs{\si_{0,1} - \si'_{0,1}},2\pi - \abs{\si_{0,1} - \si'_{0,1}}, b_1\},
 \end{gather*}
 where $(S)^\infty, (S')^\infty \in\rxyz$ with
 \begin{gather*}
 (S)^\infty = \sum_{i,j\geq0}\si_{i,j}X^iY^j\qquad \text{and} \qquad (S')^\infty = \sum_{i,j\geq0}\si'_{i,j}X^iY^j.
 \end{gather*}
 \end{Definition}

 In order for a sequence of series to converge with respect to this metric each term except for the coef\/f\/icient to $X^0Y^1$ must converge in~$\R$, and the coef\/f\/icient to $X^0Y^1$ must converge in~$\R/2\pi\Z$.

 \subsubsection{The metric on the polygon invariant}\label{sec_metricpoly}

In~\cite{PePRS2013} the authors use the Lebesgue measure of the symmetric dif\/ference of the moment polytopes to def\/ine a metric on the space of toric systems. For the portion of the semitoric metric related to the polygon invariant something similar is done in~\cite{PaSTMetric2015}.

 \begin{Definition}\label{def_admissiblemeas} Let $\nu$ be a measure on $\R^2$. We say that it is \emph{admissible} if:
 \begin{enumerate}\itemsep=0pt
 \item[1)] it is in the same measure class as the Lebesgue measure on $\R^2$ (that is, $\mu\ll\nu$ and $\nu\ll\mu$ where $\mu$ is the Lebesgue measure);
 \item[2)] 
 the Radon--Nikodym derivative of $\nu$ with respect to Lebesgue measure depends only on the $x$-coordinate, that is, there exists a $g\colon \R\to\R$ such that $\nicefrac{\text{d}\nu}{\text{d}\mu}(x,y) = g(x)$ for all $(x,y)\in\R^2$;
 \item[3)] 
 the function $g$ from part (2) satisf\/ies $xg\in\text{L}^1(\mu, \R)$ and $g$ is bounded and bounded away from zero on any compact interval.
 \end{enumerate}
 \end{Definition}

 Since the polygons may be noncompact an admissible measure on $\R^2$ is used in place of the Lebesgue measure. Notice that if~$\nu$ is an admissible measure and $\De$ is a convex noncompact polygon with everywhere f\/inite height then $\nu(\De)<\infty$. This is essentially because the convex polygon can grow at most linearly and the admissible measure is designed to shrink faster than this by part~(3). Since there is a family of polygons they must all be compared. Let $\symdiff$ denote the symmetric dif\/ference of two sets. That is, for $A,B\subset\R^2$ we have \begin{gather*} A\symdiff B = (A \setminus B)\cup(B \setminus A).\end{gather*}

 \begin{Definition} Let $\stpolymfk\subset\stpoly$ denote the set whose elements are orbits under $G_\mf\times\mathcal{G}$ of primitive semitoric polygons with twisting index $\vec{k}$ in $\Z^\mf$. Then
 \begin{gather*}\stpoly = \bigcup_{\substack{\mf\in\Z_{\geq0}\\\vec{k}\in\Z^\mf}}\stpolymfk.\end{gather*}
 \end{Definition}

 \begin{Remark} Let $\mf$ be a nonnegative integer and $\vec{k}\in\Z^\mf$. Suppose that a primitive semitoric polygon $\De_w$ has twisting
 index $\vec{k'}\in\Z^\mf$ where there exists some $c\in\Z$ such that $k_j = k'_j + c$ for $j=1, \ldots, \mf$. Then $[\De_w]\in\stpolymfk$ because the set $[\De_w]$ is also the orbit of a primitive semitoric polygon with twisting index $\vec{k}$.
 \end{Remark}

 \begin{Definition}\label{def_distpoly}
 Let $\mf\in\Z_{\geq0}$ and $\vec{k}\in\Z^\mf$. Let
 \begin{gather*}\De_w=\pstpoly \qquad \text{and}\qquad \De'_w=\big(\De', (\ell_{\la'_j}, +1, k_j)_{j=1}^{\mf}\big)\in\lwpoly\end{gather*} be primitive semitoric polygons so $[\De_w], [\De'_w]\in \stpolymfk$ are semitoric polygons in the same twisting index class. Then, if $\mf>0$ we def\/ine the distance between them to be
 \begin{gather*}
 \distpoly([\De_w],[\De'_w]) = \sum_{\vec{u}\in\{0,1\}^\mf} \nu \big(t^{\vec{u}}_{\vec{\la}}(\De) \symdiff t^{\vec{u}}_{\vec{\la'}}(\De')\big).
 \end{gather*}
 If $\mf=0$ then this sum will be empty so we instead use
 \begin{gather*}
 \distpoly([\De_w],[\De'_w]) = \inf_{k\in\Z}\big\{\nu(\De \symdiff T^k(\De'))\big\}.
 \end{gather*}
 For any $\mf\in\Z_{\geq 0}$ and for any $\vec{k}\in\Z^\mf$, the function $\distpoly\colon \stpolymfk\times\stpolymfk\to\R$ is a~metric~\cite{PaSTMetric2015}; let $\tau_\mathcal{D}$ denote the topology that it induces on $\stpolymfk$.
 \end{Definition}

The distance def\/ined in Def\/inition~\ref{def_distpoly} does not depend on the choice of primitive semitoric polygon by~\cite[Proposition~4.5]{PaSTMetric2015}, because we have specif\/ically chosen representatives which have the same twisting index~$\vec{k}$. To f\/ind the distance between two families of polygons we take the sum of the symmetric dif\/ferences of all of them and in the case that $\mf=0$ there in no twisting index $\vec{k}$ to align so we take whichever alignment is best\footnote{This is slightly dif\/ferent from the def\/inition in~\cite{PaSTMetric2015} to correct for a small mistake (which doesn't af\/fect the main results of~\cite{PaSTMetric2015}).}.

 \subsubsection{The metric on the moduli space of semitoric systems}\label{sec_fullmetric}

The metric on semitoric systems will be formed by combining the metrics from Section~\ref{sec_metricts} and Section~\ref{sec_metricpoly}.

 \begin{Definition} Let $\mmfk\subset\mmf$ denote those lists of ingredients with the polygon invariant in the set $\stpolymfk$.
 \end{Definition}

We will only def\/ine the metric on each separate $\tmfk$ because that is all we will need for this paper, see Remark~\ref{rmk_metricff}.

 \begin{Definition}[{\cite{PaSTMetric2015}}]\label{def_metric} Let $\nu$ be an admissible measure on $\R^2$ and let $\bn$ be a linearly summable sequence. Let $\mf$ be a~nonnegative integer and $\vec{k}\in\Z^\mf$.
 \begin{enumerate}\itemsep=0pt
 \item We def\/ine \emph{the metric on $\mmfk$ relative to $\nu$ and $\bn$} to be given by
 \begin{gather*} \distmmfk (m, m') = \distpoly([\De_w],[\De'_w]) + \sum_{j=1}^\mf \big( \dts((S_j)^\infty,(S'_j)^\infty) + \abs{h_j - h'_j}\big),\end{gather*}
 where $m,m'\in\mmfk$ are given by \begin{gather*}m=\melement, \qquad m'=\melementprime.\end{gather*}
 \item
 Let $\Phi\colon \mathcal{T}\to\mathcal{M}$ be the correspondence from Theorem~\ref{thm_class}, and let
 \begin{gather*}
 \mathcal{T}_{\mf,\vec{k}} = \Phi^{-1}(\mmfk).
 \end{gather*}
 The \emph{metric on $\mathcal{T}_{\mf, \vec{k}}$} is given by $\distt_{\mf,\vec{k}} = \Phi^* \distmmfk$.
 \end{enumerate}
 \end{Definition}

 \begin{Remark}\label{rmk_metricff}The metrics on each $\tmfk$ can be combined to form a metric on the whole space~$\mathcal{T}$ such that each
$\tmfk$ is an entire connected component (so the dif\/ferent $\tmfk$ either coincide or are in separate components) in a~number of ways, see~\cite{PaSTMetric2015}. For instance, we could take the distance between any two elements of $\tmfk$ to be $\min\big\{1,\distt_{\mf,\vec{k}}\big\}$ and the distance between any two elements not in the same~$\tmfk$ to be~1. For the purposes of this paper we are only interested in the topology of~$\mathcal{T}$ so we will not def\/ine a metric on all of~$\mathcal{T}$. Though it is extended to all of~$\mathcal{T}$ it is best to think of the metric from~\cite{PaSTMetric2015} as a metric and topology on each space $\tmfk$ and not the total space.
 \end{Remark}

\begin{Remark}It is important to notice that the metric we are using in this article is not the same as the metric def\/ined in~\cite{PaSTMetric2015}. Since these two metrics induce the same topology~\cite[Section~4.6]{PaSTMetric2015} it is suggested in~\cite[Remark~3.15(4)]{PaSTMetric2015} that the metric in the present article be used when studying the topological properties, such as connectedness.

The metric def\/ined in~\cite{PaSTMetric2015} produces the appropriate metric space structure, which can be seen when the completion is computed in that article. That metric is def\/ined as the minimum of a~collection of functions, one of which is the metric used in this article, so it is immediate that the distance between two systems using Def\/inition~\ref{def_metric} will never be smaller than the distance between those two systems using the metric studied in~\cite{PaSTMetric2015}.
\end{Remark}

The metric on each component $\tmfk$ depends on the choice of admissible measure and linearly summable sequence, but the topology it induces does not.

\begin{Proposition}[{\cite[Theorem A]{PaSTMetric2015}}]\label{prop_STMetricThmA}
Let $\nu$ be an admissible measure on $\R^2$ and $\bn$ a~linearly summable sequence. Then the space $\big(\mathcal{T}_{\mf,\vec{k}}, \distt_{\mf,\vec{k}}\big)$ is a metric space for any choice of $\mf\in\Z_{\geq 0}$ and $\vec{k}\in\Z^{\mf}$. Furthermore, the topology induced by
$\distt_{\mf,\vec{k}}$ on $\mathcal{T}_{\mf,\vec{k}}$ does not depend on the choice of $\nu$ or the choice of $\bn$.
\end{Proposition}

It follows from Proposition~\ref{prop_STMetricThmA} that the following notion is well def\/ined. Observe that for each $\mf\in\Z_{\geq 0}, \vec{k} \in\Z^\mf$ there exists a unique $\vec{k}'\in\Z^\mf$ with $k_1' = 0$ such that $\mathcal{T}_{\mf,\vec{k}}=\mathcal{T}_{\mf,\vec{k}'}$. Moreover, if $\mathcal{T}_{\mf,\vec{k}'} = \mathcal{T}_{\mf,\vec{k}''}$ and $k_1'=k_1''=0$ then $k'=k''$, so the union $\cup_{\vec{k}\in\Z^{\mf},k_1 = 0}(\mathcal{T}_{\mf,\vec{k}})$ is the set of isomorphism classes of semitoric systems with $\mf$ focus-focus points and any two distinct sets in the above union are disjoint.

\begin{Definition}
The \emph{topology on $\mathcal{T}$} is the disjoint union topology on
\begin{gather*}
 \mathcal{T} = \bigsqcup_{\substack{\mf\in\Z_{\geq 0},\\ \vec{k}\in\Z^{\mf}\textrm{ with }k_1 = 0}}\mathcal{T}_{\mf,\vec{k}},
\end{gather*}
where the topology on each $\tmfk$ is induced by the metric $\distt_{\mf,\vec{k}}$ for a choice of admissible metric $\nu$ and linear summable sequence~$\bn$.
\end{Definition}

This topology is def\/ined so that each $\tmfk$ is a separate component of $\mathcal{T}$, so it is natural to wonder if $\tmfk$ is connected for each f\/ixed $\mf\in\Z_{\geq0}$ and $\vec{k}\in\Z^\mf$. The remainder of the paper is devoted to proving Theorem~\ref{thm_stconnect}, which states that these are in fact path-connected.

\subsection{The connectivity of the moduli space of semitoric integrable systems}\label{sec_stconnect}

\begin{Definition}\label{def_assocstfan}
 Let $\De_w\in\lwpoly$ be a compact primitive semitoric polygon with
 \begin{gather*}
 \De_w = \pstpoly.
 \end{gather*}
 Then the \emph{associated semitoric fan} is the semitoric fan $\mathcal{F}$ formed by the inwards pointing primitive integer normal vectors to the edges of $\De$ in which the pair of vectors in $\mathcal{F}$ are labeled as fake, hidden, or Delzant to correspond with the labeling of the corners of $\De$.
\end{Definition}

Let $\nu$ be an admissible measure (cf.\ Def\/inition~\ref{def_admissiblemeas}), endow the set $\poly$ with the topology induced by the metric given by the $\nu$-measure of the symmetric dif\/ference, and endow $\poly\times\R^\mf\times\Z^\mf$ with the product topology. The topology induced on $\poly$ does not depend on the choice of admissible measure and the topology induced on the set of compact rational convex polygons by the Lebesgue measure of the symmetric dif\/ference agrees with the subset topology induced by the topology on~$\poly$, see~\cite{PaSTMetric2015}.
\begin{Definition}\label{def_topologyprimitive}
The set of primitive semitoric polygons of defect $\mf\in\Z_{\geq0}$ inherits a topology from the bijection $\big(\De, (\ell_{\la_j}, +1, k_j)_{j=1}^{\mf}\big) \mapsto \big(\De, (\la_j)_{j=1}^\mf, (k_j)_{j=1}^\mf\big)$, which we denote by $\tau_0$.
\end{Definition}

\begin{Lemma}\label{lem_continuoustransf} Every semitoric fan is the associated fan of some compact primitive semitoric polygon, and each relation in Theorem~{\rm \ref{thm_stpoly}} corresponds to some continuous transformation of compact primitive semitoric polygons. More specifically, suppose that two semitoric fans $\mathcal{F}_0, \mathcal{F}_1 \in \big(\Z^2\big)^d$ are related by
 \begin{enumerate} \itemsep=0pt
 \item[$1)$] performing corner chops;
 \item[$2)$] performing reverse corner chops;
 \item[$3)$] removing hidden corners; or
 \item[$4)$] commuting fake and Delzant corners;
 \end{enumerate}
$($see Definition {\rm \ref{def_fantrans})}. Then there exists a continuous $($with respect to~$\tau_0$, Definition~{\rm \ref{def_topologyprimitive})} family of (compact) primitive semitoric polygons
 \begin{gather*}\De_w^t=\big(\De_t, (\ell_{\la^t_j}, +1, k_j)_{j=1}^{\mf}\big),\end{gather*}
 $t\in[0,1]$ with fixed $k_1, \ldots, k_\mf$, such that the semitoric fan associated to $\De_w^0$ is $\mathcal{F}_0$ and the fan associated to $\De_w^1$ is $\mathcal{F}_1$.
\end{Lemma}

\begin{proof}Given any semitoric fan it is possible to construct a polygon $\De$ which has those vectors as inwards pointing normal vectors, and by the def\/inition of semitoric fans $(\De, (\ell_{\la_j}, +1, 0)_{j=1}^c)$ (where $\la_j$ is the $x$-value of the location of the $j^{th}$ hidden or fake corner) is a primitive semitoric polygon with the prescribed associated semitoric fan. Suppose that $(\De, (\ell_{\la_j}, +1, k_j)_{j=1}^{\mf})$ is a~compact primitive semitoric polygon with associated semitoric fan $\mathcal{F} = (v_0, \ldots, v_{d-1})\in\big(\Z^2\big)^d$
and f\/ix some $i\in\{0, \ldots, d-1\}$. Let $v_{-1}:=v_{d-1}$ and $v_d:=v_{0}$ so that the formulas used in this proof will be valid if $i=0$ or $i=d-1$. Throughout, let $p\in\De$ be the corner adjacent to the edges with inwards pointing normal vectors $v_i$ and $v_{i+1}$. Let $u_i\in\Z^2$, $i=1,2$, denote the primitive vectors along which the edges adjacent to $p$ are aligned (oriented to point outwards from $p$), ordered so that
$\det(u_1,u_2)>0$.

For $w_1, w_2\in\Z^2$ let $\mathcal{H}_p^\varep (w_1,w_2)$ denote the half-space given by
\begin{gather*}
 \mathcal{H}_p^\varep (w_1,w_2) = \{ p+t_1 w_1 + t_2 w_2\colon t_1 +t_2 \geq \varep \}.
\end{gather*}

First we consider the corner chop operation. Suppose $p$ is a Delzant corner. Fix some $\varep_0>0$ smaller than the length of the edges incident at $p$.

For $t\in[0,1]$ let
\begin{gather*}
 \De_t = \De \cap \mathcal{H}_p^{t \varep_0} (u_1,u_2).
\end{gather*}
We see that
$\De_t$ is a continuous family and since the edges of $\De_t$ are parallel to the edges of $\De$ except for the new edge with inwards pointing normal vector given by $v_i + v_{i+1}$ we see that the semitoric fan of $\De_t$ is the corner chop of the semitoric fan for $\De$ for $t\in(0,1]$. Since a reverse corner chop is the inverse of this operation, we can use the same path backwards.

Now suppose that $p$ is a hidden corner. Let $\be\in\Z$ be the second component of $v_{i+1}$. For $\varep_0>0$ smaller than the length of the adjacent edges and
$t\in[0,1]$ let
\begin{gather*}
 \De_t = \De \cap \mathcal{H}_p^{t\varep_0} \big(u_1, \be^2 u_2\big).
\end{gather*}
The integral normal vector to the new edge of the polygon is $Tv_{i+1}$. Indeed let $v_i = \left(\begin{smallmatrix} a\\ b\end{smallmatrix}\right)$ and $v_{i+1}=\left(\begin{smallmatrix} \al \\ \be\end{smallmatrix}\right)$, so the new vector should be parallel to
\begin{gather*}
\be^2 v_i + v_{i+1} = \be^2 v_i + (a\be - \al b - b\be) v_{i+1} =\vect{\be^2 a+\al (a\be - \al b - b\be)}{\be^2 b + \be (a\be - \al b - b\be)} = (a\be - \al b) T v_{i+1}
\end{gather*}
using the fact that $\det(v_i,Tv_{i+1}) = a\be - \al b - b\be = 1$ since $p$ is a hidden corner. We see that~$\De_t$ is a continuous family and by construction it has the desired semitoric fan. Thus $\De_t$ is the required family for the operation of removing hidden corners.

Finally, suppose that $p$ is a fake corner and the next corner, which has adjacent edges which have inwards pointing normal vectors $v_{i+1}$ and $ v_{i+2}$, is a Delzant corner on the top boundary. Since $v_i$, $v_{i+1}$, and $v_{i+2}$ are all in the lower half plane we see $d>3$ since a semitoric fan must always include at least one vector in the upper half plane. Viewing the convex polygon as an intersection of half-planes, remove the half-plane corresponding to the edge with inwards pointing normal vector $v_{i+1}$ from the intersection to produce a new polygon $\De'$ which has one fewer edge. That is, extend the edges with inwards pointing normal vectors $v_i$ and $v_{i+2}$ until they intersect; call the new corner where they intersect $p'$. The polygon $\De'$ has semitoric fan
\begin{gather*}
 \mathcal{F}' = (v_0, \ldots, v_i, v_{i+2}, \ldots, v_{d-1})\in\big(\Z^2\big)^{d-1}.
\end{gather*}
Since $(v_i, v_{i+1})$ is fake and $(v_{i+1}, v_{i+2})$ is Delzant by assumption we have
\begin{gather*}
 \det(v_i, Tv_{i+2}) = \det(Tv_{i+1}, Tv_{i+2}) = \det(v_{i+1}, v_{i+2})=1,
\end{gather*}
so $\mathcal{F}'$ is a semitoric fan and the corner $(v_i, v_{i+2})$ in $\mathcal{F}'$ is a hidden corner. Removing the hidden corner $(v_i, v_{i+2})$ from $\mathcal{F}'$ produces the required semitoric fan $(v_0, \ldots, v_i, Tv_{i+2}, v_{i+2},\ldots, v_{d-1})$, and since removing a hidden corner can be achieved by a continuous family (as shown above) all that remains is to produce a continuous family from a polygon with semitoric fan $\mathcal{F}$ to~$\De'$. Let~$u_1'$ be the primitive integral vector directing the edge away from $p'$ which has inwards pointing
normal vector $v_{i+2}$ and let $u_2'$ be the primitive integral vector directing the other edge away from~$p'$. Then for some suf\/f\/iciently small $\varep_0>0$ let
\begin{gather*}
 \De_t = \De' \cap \mathcal{H}_{p'}^{t\varep_0}(u_1', u_2'),
\end{gather*}
and note that $\De_0 = \De'$ and $\De_1$ has the same semitoric fan as $\De$ (for the correct choice of $\varep_0$ we can actually get $\De_0 = \De$, but this is not necessary), so we have formed a~continuous path from a polygon with associated semitoric fan $\mathcal{F}$ to the polygon $\De'$, as desired.

So in any case we have produced a family of polygons $\De_t$, $0\leq t \leq 1$. Let
\begin{gather*}
 \De_w^t = \big(\De_t, (\ell_{\la_j^t}, +1, k_j)_{j=1}^\mf\big),
\end{gather*}
where $\la_j^t$ is the $x$-coordinate of the $j^{\textrm{th}}$ non-Delzant corner of $\De_t$ ordered left to right, to complete the proof.
\end{proof}

\begin{figure}\centering
 \includegraphics[width=300pt]{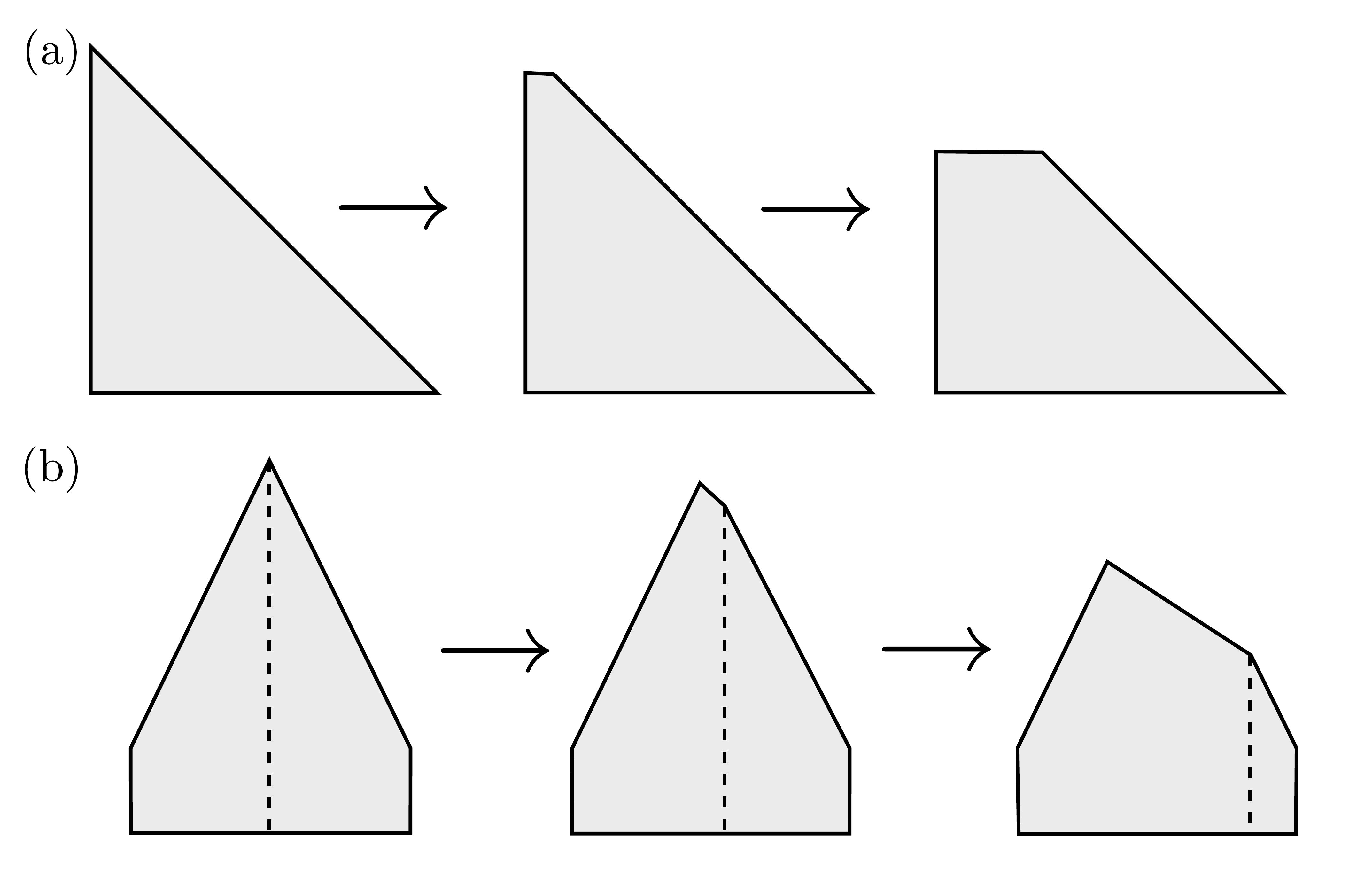}
 \caption{For each fan transformation there is a continuous path of semitoric polygons which transitions between the semitoric fans. In (a) we show
 the corner chop and in (b) we show the removal of a hidden corner (which is replaced by a fake and a Delzant corner).} \label{fig_contchops}
\end{figure}

Recall that toric polygons are precisely the compact primitive semitoric polygons with defect zero. Thus, by Proposition~\ref{prop_standardformstfan} and Lemma~\ref{lem_continuoustransf} we have recovered Theorem~\ref{thm_toricconnected}.

In light of Propositions~\ref{prop_standardformstfan} and~\ref{lem_continuoustransf} the only dif\/f\/iculty remaining to prove the following lemma is incorporating the case of semitoric systems which have noncompact polygons as invariants. Recall $\tau_\mathcal{D}$ is the topology on $\stpolymfk$ induced by the metric def\/ined in~\ref{def_distpoly}.

\begin{Lemma}\label{lem_dpolypathconn} Let $\mf\in\Z_{\geq0}$ and $\vec{k}\in\Z^\mf$. Then $\stpolymfk$ is path-connected with respect to the topology $\tau_\mathcal{D}$.
\end{Lemma}

\begin{proof} Every semitoric polygon is an orbit of a primitive semitoric polygon, and to conclude that there is a continuous path between two semitoric polygons it is suf\/f\/icient to construct a~continuous (with respect to $\tau_0$, Def\/inition~\ref{def_topologyprimitive}) path between the primitive semitoric polygons that generate those orbits. This is because the map which takes a primitive semitoric polygon to its orbit is continuous. Thus, it suf\/f\/ices to consider only primitive semitoric polygons. Any two compact primitive semitoric polygons in $\stpolymfk$ with the same semitoric fan can be connected by a continuous path. This path is made by continuously changing the lengths of the edges because the angles of the two polygons must all be the same since they have the same semitoric fan while also changing $\la_j$ to still correspond to the locations of the fake and hidden corners. So, by Proposition \ref{prop_standardformstfan} given two elements of $\stpolymfk$ which are compact we know that the corresponding semitoric fans are related by the moves listed in that proposition and then by Lemma~\ref{lem_continuoustransf} we know these moves correspond to continuous paths of polygons. So we have established that any two compact elements of $\stpolymfk$ are connected by a~continuous path.

 Next assume that
 \begin{gather*}
 [\De_w] = \big[\lwp\big]\in\stpolymfk
 \end{gather*}
 is such that $\De$ is noncompact but has only f\/initely many vertices. Let
 \begin{gather*}
 \pi_1\colon \ \R^2\to\R
 \end{gather*}
 be projection onto the f\/irst component. Choose $N\in\R$ such that all of the vertices of $\De$ are in the region $\pi_1^{-1}([-N,N])$. The set $\De \cap \pi_1^{-1}([-N-1, N+1])$ is a polygon but the corners which intersect $\ell_{N+1}\cup\ell_{-N-1}$ may not satisfy the Delzant condition.
 By~\cite[Remark 23]{PePRS2013} we may change the set on arbitrarily small neighborhoods of these corners to produce a new compact polygon, $\De'$, which is equal to $\De \cap \pi_1^{-1}([-N-1, N+1])$ outside of those small neighborhoods and has only Delzant corners inside of those neighborhoods. Thus,
 \begin{gather*}
 \De'_w = \big(\De', (\ell_{\la_j}, \ep_j, k_j)_{j=1}^\mf\big)
 \end{gather*}
 is a primitive semitoric polygon and by choosing the neighborhoods small enough we can assure that $\De \cap \pi_1^{-1}([-N,N]) = \De' \cap \pi_1^{-1}([-N,N])$. For $t\in[0,1)$ let $\De(t)$ be the polygon with the same semitoric fan as~$\De'$, the property that
 \begin{gather*}
 \De(t) \cap [-N,N] = \De' \cap [-N,N],
 \end{gather*}
 and which has all of the same edge lengths as $\De'$ with the exception of the two or four edges which intersect $\ell_N\cup\ell_{-N}$. These
 edges are extended horizontally by a length of $\frac{1}{t-1}$. By this we mean that if an edge of $\De'$ which intersected $\pi_1^{-1}(N)$ had as one of its endpoints $(x_0, y_0)$ with $x_0>N$ then the corresponding edge of $\De(t)$ would have as its endpoint $\big(x_0+\frac{1}{t-1}, y_0+m\frac{1}{t-1}\big)$, where $m$ is the slope of the edge in question. Then def\/ine $\De(1) = \De$ and we can see that $\De(t)$ for $t\in[0,1]$ is a path from $\De'$, which is compact, to $\De$ so $\big[\big(\De(t), (\ell_{\la_j}, +1, k_j)_{j=1}^\mf\big)\big]$ is a~continuous path which connects a compact semitoric polygon to $[\De_w]$. This process is shown in Fig.~\ref{fig_dpolypath1}.
\begin{figure}[t]\centering
 \includegraphics[width=280pt]{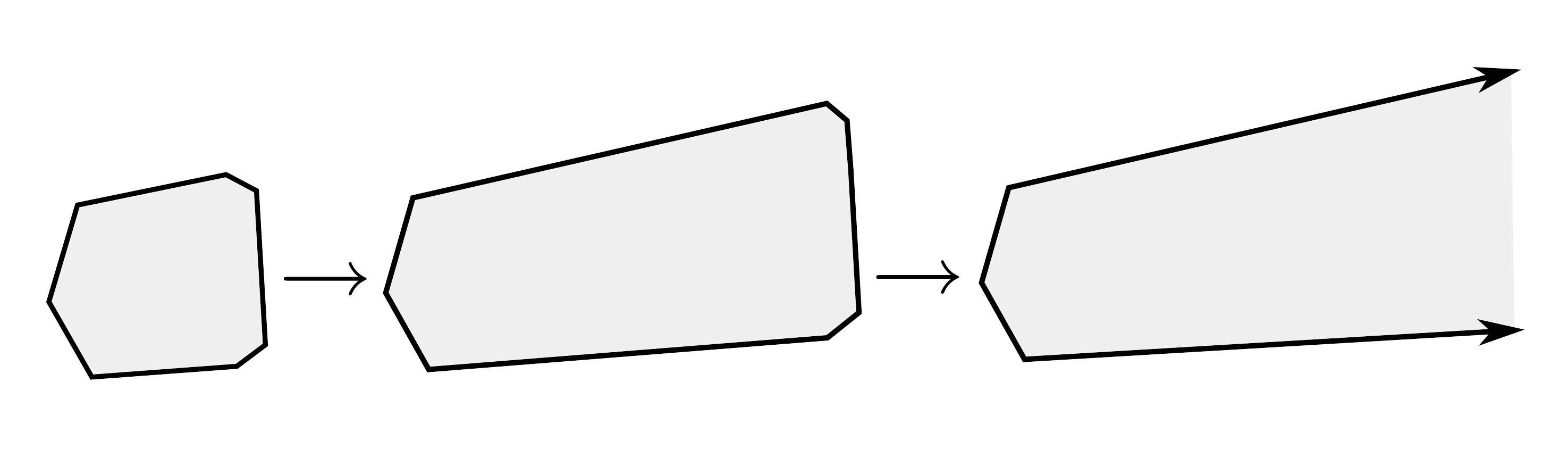}
\caption{The continuous path from a compact primitive semitoric polygon to a noncompact primitive semitoric polygon with f\/initely many vertices.}
 \label{fig_dpolypath1}
\end{figure}

Now we have connected all of the elements except for those with an inf\/inite amount of vertices. Suppose that
 \begin{gather*} [\De_w] = \big[\lwp\big]\in\stpolymfk\end{gather*}
is such that $\De$ is noncompact and has inf\/initely many vertices. We will connect $[\De_w]$ to a semitoric polygon which has only f\/initely many vertices to f\/inish the proof. Since $\De$ has everywhere f\/inite height and is the intersection of inf\/initely many half-planes we can choose two of these planes which are not horizontal and are not parallel to one another. Denote the intersection of these two half-planes by $A$ and notice $\De\subset A$. Since the boundaries of these two half-planes must intersect we can see that $A$ can only be unbounded in either the positive or negative $x$-direction, but not both. Without loss of generality assume that $A$ is unbounded in the positive $x$-direction.

Let $\nu$ be any admissible measure. For any $n\in\Z_{\geq0}$ since $\nu(A)<\infty$ we know there exists some $x_n \in \R$ such that
\begin{gather*}
 \nu\big(A\cap\pi_1^{-1}([x_n, \infty))\big)<\nicefrac{1}{n},
\end{gather*}
$\De$ does not have a corner on the line $\ell_{x_n}$, and $x_n>\abs{\la_j}$ for all $j=1, \ldots, \mf$. Let $\De_n$ denote the polygon which satisf\/ies
\begin{gather*}
 \De_n \cap \pi_1^{-1}([-\infty, x_n]) = \De \cap \pi_1^{-1}([-\infty,x_n])
\end{gather*}
and has no vertices with $x$-coordinate greater than $x_n$.

\begin{figure}[t]\centering
 \includegraphics[width=280pt]{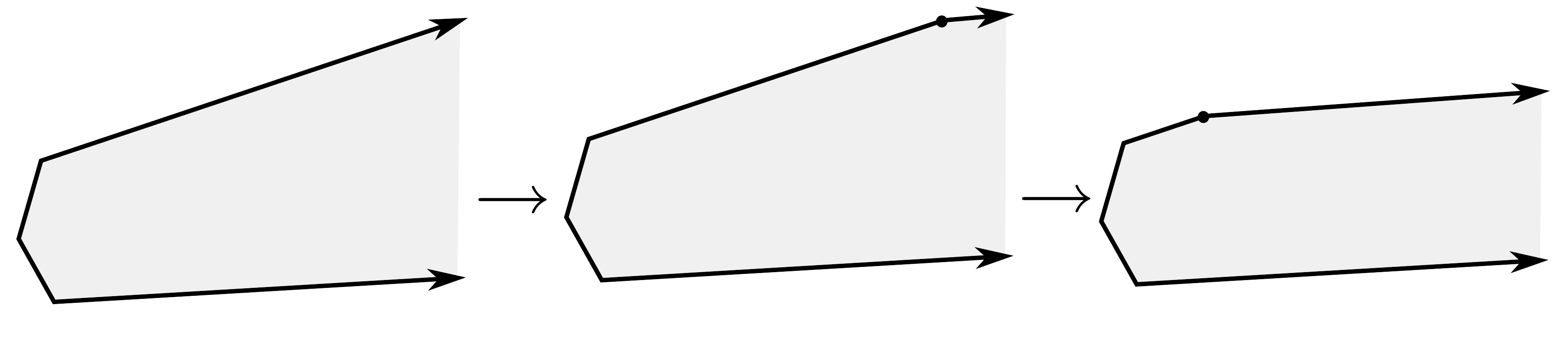}
\caption{The continuous path which adds new vertices to a~noncompact primitive semitoric polygon. This process is repeated to form a path to systems with inf\/initely many vertices.} \label{fig_dpolypath2}
\end{figure}

For each $n\in\Z_{\geq0}$ and $t\in(0,1]$ def\/ine $\De_n(t)$ to have the same semitoric fan as $\De_{n+1}$ and to have all the same edge lengths as $\De_{n+1}$ except for the two edges which intersect $\ell_{x_n}$. Extend those two edges horizontally by $\nicefrac{1}{t}-1$. Def\/ine $\De_n(0) = \De_n$. Now $\De_n(t)$ for $t\in[0,1]$ is a $C^0$ path, with respect to the topology on polygons generated by the $\nu$-measure of the symmetric
dif\/ference, which takes $\De_n$ to $\De_{n+1}$. Moreover,
\begin{gather*}
 \De\symdiff\De_n(t) \subset A\cap\pi_1^{-1}([x_n, \infty])\qquad \text{so} \qquad \nu(\De\symdiff\De_n(t))<\nicefrac{1}{n}
\end{gather*}
for each $t\in[0,1]$. Each of these paths for $n\in\Z_{\geq0}$ can be concatenated to form a continuous path $\De_t$, $t\in[0,1]$, from $\De_0$ to $\De$ and we know that $\De_0$ has only f\/initely many vertices. It is important not only that each $\De_n$ be getting closer to $\De$ but also that the path from~$\De_n$ to~$\De_{n+1}$ stays close to~$\De$. Then we def\/ine $[(\De(t), (\la_j, +1, k_j)_{j=1}^\mf)]$ which is a continuous
path of semitoric polygons from a semitoric polygon with f\/initely many vertices to $[\De_w]$. This is shown in Fig.~\ref{fig_dpolypath2}.
\end{proof}

Now we can classify the connected components of $\m$ and $\mathcal{T}$. Recall that $\mmfk = \mathcal{M}_{\mf, \vec{k'}}$ if $k_j = k'_j +c$ for some $c\in\Z$ so when stating the following lemma we require that the f\/irst component of the twisting index be $0$. This is done only to make
sure that there are no repeats in the list of components. Recall $\m_0$ is the collection of semitoric lists of ingredients with $\mf=0$ and $\mathcal{T}_0 = \Phi^{-1}(\m_0)$ is the collection of semitoric systems with no focus-focus singularities. Finally, also recall that the fact that the dif\/ferent $\mmfk$ are in dif\/ferent components really comes from the def\/inition of the constructed topology in~\cite{PaSTMetric2015} and not from any result of the present paper.

\begin{Lemma}\label{lem_stpolyconnect}
 The connected components of $\m$ are each element of the set
 \begin{gather*}\{\mmfk\,|\,\mf\in\Z_{>0},\vec{k}\in\Z^\mf\textrm{ with }k_1 = 0\}\cup\{\m_0\}\end{gather*} and they are each path-connected.
\end{Lemma}

\begin{proof} It is suf\/f\/icient to prove that $\mmfk$ is path-connected for each choice of $\mf\in\Z_{\geq0}$ and $\vec{k}\in\Z^\mf$. Let $m,m'\in\mmfk$ with \begin{gather*}m = \melement \qquad \text{and} \qquad m'=\melementprime. \end{gather*} By Lemma \ref{lem_dpolypathconn} we know there exists a continuous path \begin{gather*}[\De_w(t)]=\big[\big(\De(t), (\ell_{\la_j(t)}, +1, k_j)_{j=1}^\mf\big)\big],\end{gather*} $t\in[0,1]$, from $[\De_w]$ to $[\De'_w]$ and by~\cite[Proposition 3.2]{PaSTMetric2015} we know $\rxyz$ is path-connected so there exists a continuous path $(S_j(t))^\infty$ from $(S_j)^\infty$ to $(S'_j)^\infty$ for each $j=1, \ldots, \mf$. For $j=1, \ldots, \mf$ let
 \begin{gather*}
 \textrm{len}_j = \textrm{length}(\pi_2 (\De \cap \ell_{\la_j })),\qquad
 \textrm{len}'_j = \textrm{length}(\pi_2 (\De' \cap \ell_{\la'_j })),\\
 \textrm{len}_j(t) = \textrm{length}(\pi_2 (\De(t) \cap \ell_{\la_j(t)}))
 \end{gather*}
 for $t\in[0,1]$ and def\/ine
 \begin{gather*}
 h_j(t) = \left(\frac{(1-t)h_j}{\textrm{len}_j} + \frac{t h_j'}{\textrm{len}'_j}\right)\textrm{len}_j(t).
 \end{gather*}
 Now we have that $0<h_j(t)<\textrm{len}_j(t)$ and $t\mapsto h_j(t)$ is a continuous function from $[0,1]$ to~$\R$ because it is impossible for a semitoric polygon to have a vertical boundary at $\ell_{\la_j}$ for any $j\in\{1, \ldots, \mf\}$. Now def\/ine \begin{gather*}m(t) = \big([\De_w(t)], (h_j(t))_{j=1}^{\mf}, ((S_j(t))^\infty)_{j=1}^{\mf}\big) \end{gather*} for $t\in[0,1]$ which is a continuous path from~$m$ to~$m'$.
\end{proof}

Thus we have established the following result.

\begin{Theorem}\label{thm_stconnect2} The set of connected components of $\mathcal{T}$ is
 \begin{gather*}
 \big\{\tmfk\,|\,\mf\in\Z_{>0},\, \vec{k}\in\Z^\mf\textrm{ with }k_1 = 0\big\}\cup \{\mathcal{T}_0\}
 \end{gather*}
 and they are each path-connected.
\end{Theorem}

Theorem \ref{thm_stconnect} is equivalent to Theorem \ref{thm_stconnect2}.

\subsection*{Acknowledgements}

We thank the anonymous referees for reading the paper carefully and providing very helpful suggestions which have improved the paper.
JP and \'AP were partially supported by NSF grants DMS-1055897 and DMS-1518420.

\pdfbookmark[1]{References}{ref}
\LastPageEnding

\end{document}